\documentclass[preprint,review,3p,fleqn]{elsarticle}

\usepackage[linesnumbered,ruled]{algorithm2e}
\usepackage{graphics,graphicx,epsfig}
\usepackage{amsmath,amsfonts,amssymb}
\usepackage{xcolor}
\usepackage{color}
\usepackage{framed}
\usepackage{bm}
\usepackage{graphicx}
\usepackage{booktabs,multirow}
\usepackage{caption}
\usepackage{lineno}
\let\origlinenumberfont\linenumberfont
\renewcommand{\linenumberfont}{\origlinenumberfont}

\newcommand{\be}{\begin{eqnarray*}}
\newcommand{\ee}{\end{eqnarray*}}
\newcommand{\ibe}{\begin{eqnarray}}
\newcommand{\iee}{\end{eqnarray}}

\allowdisplaybreaks

\begin{document}

\begin{frontmatter}

\title{Three-dimensional narrow volume reconstruction method with unconditional stability based on a phase-field Lagrange multiplier approach}
\author[math1]{Renjun Gao}
\author[math1]{Xiangjie Kong}
\author[math1]{Dongting Cai}
\author[math1]{Boyi Fu}
\author[math1]{Junxiang Yang\corref{cor1}}\ead{jxyang@must.edu.mo}
\cortext[cor1]{Corresponding author.}
\address[math1]{School of Computer Science and Engineering, Faculty of Innovation Engineering, Macau University of Science and Technology, Macao Special Administrative Region of China}

\begin{abstract}
Reconstruction of an object from points cloud is essential in prosthetics, medical imaging, computer vision, etc. We present an effective algorithm for an Allen--Cahn-type model of reconstruction, employing the Lagrange multiplier approach. Utilizing scattered data points from an object, we reconstruct a narrow shell by solving the governing equation enhanced with an edge detection function derived from the unsigned distance function. The specifically designed edge detection function ensures the energy stability. By reformulating the governing equation through the Lagrange multiplier technique and implementing a Crank--Nicolson time discretization, we can update the solutions in a stable and decoupled manner. The spatial operations are approximated using the finite difference method, and we analytically demonstrate the unconditional stability of the fully discrete scheme. Comprehensive numerical experiments, including reconstructions of complex 3D volumes such as characters from \textit{Star Wars}, validate the algorithm's accuracy, stability, and effectiveness. Additionally, we analyze how specific parameter selections influence the level of detail and refinement in the reconstructed volumes. To facilitate the interested readers to understand our algorithm, we share the computational codes and data in https://github.com/cfdyang521/C-3PO/tree/main.
\end{abstract}


\end{frontmatter}


\section{Introduction}\label{sec1}

In both medical and industrial fields, one of the important techniques is three-dimensional (3D) printing which demonstrates considerable potential in additive manufacturing, biomedical modeling, and customized design. The reconstruction of 3D object is a crucial technique for printing \cite{shapeRe01,shapeRe02}. Since printers require precise 3D models to fabricate objects, they rely on volume reconstruction methods to convert original data into printable models. Moreover, 3D reconstruction techniques are widely applied across numerous disciplines, such as computer graphics \cite{3Dimmersive1}, medical imaging using tomographic slices \cite{MMimage0}, and computational fluid dynamics \cite{JXsurface0,JXsurface1}. For instance, in computer graphics, 3D reconstruction enables the creation of realistic virtual environments; in medical imaging, it facilitates the modeling of anatomical structures by reconstructing successive image slices; in materials science, it aids in analyzing microstructures through detailed 3D models; and in computational fluid dynamics, it assists in simulating flow patterns around complex geometries.

Point cloud technique offers a highly promising solution for 3D volume reconstruction \cite{3DpointCloud1}. Utilizing devices such as laser scanners and probe digitizers, it is possible to capture numerous data points on an object's surface, thereby accurately representing its 3D structure as a point cloud. Reconstruction from these point cloud data enables the rapid and precise restoration of the object's geometric information while effectively handling complex surface details \cite{3DpointCloud2}. Consequently, point cloud-based reconstruction methods have emerged as a prominent research focus within the field of 3D reconstruction. These techniques are extensively applied in industrial design and medical modeling, where high precision and intricate geometric details are essential.

3D reconstruction methods are broadly categorized into explicit and implicit approaches, each offering distinct techniques for handling real-world data. Explicit methods, such as Delaunay triangulation and alpha shapes \cite{Edelsbrunner1994}, represent surfaces by directly connecting data points based on proximity and local geometry. While effective with well-organized data, these methods often struggle with unstructured point clouds due to the lack of predefined connectivity, inherent noise, and non-uniform sampling. The scattered and unordered nature of points complicates the explicit description of surface geometry, posing significant challenges for complex topologies and large datasets. To mitigate these issues, Wu et al. \cite{GCN} utilized Graph Convolutional Network (GCN) to construct a diffusion model for directly extracting textured meshes from images, improving robustness even with limited input views. This approach enhances the feasibility of explicit surface reconstruction with sparse data.

The raw point cloud data often presents challenges such as disorder and noise, making accurate reconstruction difficult. To implicitly construct a smooth 3D volume, the Allen--Cahn (AC) model \cite{ACcmame01} is particularly effective due to three main properties: (i) surfaces are described by the zero level set of a characteristic function; (ii) the noise can be effectively removed via diffusion; (iii) the model's dynamics drives the initial surface closer to the points cloud.Compared with classical level–set front‑propagation frameworks \cite{Osher1988,Sethian1999,OsherFedkiw2003}, our phase‑field formulation is tied to an explicit variational energy and thus requires no signed‑distance reinitialization to enforce $|\nabla d|=1$, which reduces parameter sensitivity and geometric drift and is advantageous for noisy, unorganized point clouds \cite{Peng1999,Li2010}.
Based on the above mentioned properties, the AC phase-field theory has been widely adopted in data classification and assimilation \cite{JSKJWYB02,JSKJWYB022}, dendritic growth \cite{JSKJWYB04,JSKJWYB041}, fracture \cite{face1,face2}, crystal dynamics \cite{CCYYY01}, computational geometry \cite{acstruc2}, image inpainting and shape transformation\cite{image001,image002,image003}, etc. Over the past decade, numerous effective numerical methods for AC-type equations have been developed. For instance, Yang et al. \cite{ACe1} developed fast, unconditionally stable algorithms for the AC equation using an operator splitting method, incorporating techniques to maintain uniform interfacial transition layers under giant time steps. In recent years, numerous effective method have been proposed to address the Cahn--Hilliard, AC, and other multi-physics coupled phase-field systems; see \cite{ACe4,ACe5,ACe6,ACe7,ACe8,ACe9} for examples of such research. The scalar auxiliary variable approach provided an efficient way to design energy-stable scheme for most phase-field systems, see \cite{CCWW01,CCWW02,CCWW03} for some successful applications. However, most auxiliary variable methods modify the original energy property because an auxiliary variable replaces the nonlinear energy. To address this issue, the Lagrange multiplier approach has been explored to preserve the original energy while ensuring energy stability. Cheng et al. \cite{Cheng2020} introduced a a novel structure-preserving auxiliary variable method for phase-field models. This enhancement resulted in robust and efficient numerical schemes preserving energy stability and the original structures.

We herein develop a practical linear algorithm for generating 3D volume using the phase-field AC model, achieving second-order accuracy and unconditional energy stability through the Lagrange multiplier approach. We employ an edge detection function based on the unsigned distance function to generate narrow volume based on points cloud. By introducing a time-dependent variable $Q$ as a Lagrange multiplier, we reformulate the AC model into an equivalent system that preserves the original energy function. Specifically, we multiply $Q$ with the nonlinear term and introduce an auxiliary equation linking $Q$ to the energy function, ensuring the system is well-posed and energy-stable. This reformulation allows us to develop a unconditionally stable scheme with Crank--Nicolson (CN) method. We propose a fully discrete finite difference scheme and prove its unconditional energy stability analytically. Unlike existing methods, our approach maintains the original energy landscape without introducing auxiliary variables that alter it. The main scientific contributions of this work are as follows:\\
(i) The proposed phase-field model can effectively reconstruct various 3D volumes with complex shapes;\\
(ii) The proposed numerical solver is efficient in the sense of linear algorithm structure;\\
(iii) The proposed numerical solver is unconditionally stable in the sense of energy dissipation.

\textit{Novelty relative to prior work.}
Beyond classical pipelines (e.g., Poisson/MLS/variational families summarized in recent surveys~\cite{HuangSurvey2022,SulzerSurvey2023}), we position our contribution with respect to two representative modern directions.
First, for \emph{deterministic} reconstruction from raw or unoriented points, approaches such as variational implicit point-set surfaces and iterative Poisson for unoriented inputs~\cite{VIPSS2019,iPSR2022} substantially enhance robustness, yet they do not provide an evolution that preserves a variational energy-dissipation law.
Second, \emph{neural implicit} methods fit SDF/occupancy functions from point clouds~\cite{DeepSDF2019,OccNet2019,SAL2020,IGR2020,NeuralPull2021,Points2Surf2020}, often achieving high-fidelity reconstructions but at the expense of training data, tuning of multiple hyperparameters, and the absence of a built-in energy stability guarantee.
In contrast, we adopt an edge-weighted Allen--Cahn phase-field formulation tailored to unorganized point sets, where $g(\mathbf{x})$ is constructed from the unsigned distance to the data-driven initial interface; the original energy structure is retained so that $\tfrac{d}{dt}E(\phi) = -\lVert\sqrt{g}\,\mu\rVert^{2} \le 0$. Practically, this allows large time steps and obviates the maintenance of a signed‑distance field and ad‑hoc narrow bands that are commonly required in level‑set pipelines \cite{Sethian1999,Peng1999}.
On the algorithmic side, building upon the Lagrange-multiplier reformulation (Eqs.~\eqref{AC1}--\eqref{Qdf}), we employ an \emph{affine-in-$Q$} decoupling, namely $\phi^{n+1} = \phi^{n+1}_{1} + Q^{n+\frac12}\phi^{n+1}_{2}$ and $\mu^{n+\frac12} = \mu^{n+\frac12}_{1} + Q^{n+\frac12}\mu^{n+\frac12}_{2}$, which reduces the Crank--Nicolson update to two linear subproblems plus a single scalar equation for $Q$ determined by the discrete energy constraint (Eq.~\eqref{cnQdf}), yielding the polynomial relation (Eq.~\eqref{simplified}). Unlike SAV-type strategies, the nonlinear energy is not modified (cf.~\cite{Cheng2020}).
Practically, we obtain stable reconstructions from complex, noisy point clouds with Practically, we obtain stable reconstructions from complex, noisy point clouds with \emph{unconditional} discrete energy dissipation (Eqs.~\eqref{eq:discrete_energy_dissipation}--\eqref{eq:discrete_energy}), without any training data; we also quantify how $\epsilon$ balances geometric detail and smoothing in the recovered narrow volumes.

The following portions of this work are organized as outlined below: Section \ref{sec2} introduces the mathematical model and discuss its energy structure. In Sec. \ref{sec3}, an unconditionally stable discrete scheme utilizing the Lagrange multiplier technique is proposed. Extensive simulation will be implemented in Sec. \ref{sec4} to verify the capability of our method. Ultimately, In Sec. \ref{sec5}, we summarize our findings and conclude the study.

\section{Mathematical model}\label{sec2}
In the phase-field approach, two phases are distinguished by the characteristic function $\phi = \phi(\mathbf{x}, t)$. The temporal variable is $t$. In the bulks, $\phi$ takes values of $1$ or $-1$. Throughout the interface between two materials, $\phi$ takes the value in $(-1,1)$. In terms of $\phi$, its free energy functional is \ibe
E(\phi) = \int_\Omega \left[ \frac{F(\phi)}{\epsilon^2} + \frac{1}{2}|\nabla\phi|^2 \right] ~d{\bm x}. \label{Oene1}\iee

In this context, $\Omega$ represents the spatial domain. The nonlinear term $F(\phi) = 0.25(\phi^2-1)^2$ characterizes the separation between the two materials. A small positive constant $\epsilon > 0$ is associated with the size of the transition layer. The gradient operator is denoted by $\nabla$.
By taking the variational approach for Eq. \eqref{Oene1}, we get $\mu = \frac{F'(\phi)}{\epsilon^2} - \Delta\phi$. By applying the $L^2$-gradient flow, the AC model reads as\ibe
\frac{\partial\phi}{\partial t} &=& -\mu,\label{ac1}\\
\mu &=& \frac{F'(\phi)}{\epsilon^2} - \Delta\phi. \label{ac2}\iee
The AC equation was initially formulated to describe the anti-phase coarsening of binary alloys \cite{allencahn} by modeling the evolution of interfaces between phases.

\begin{figure}[htbp]
    \centering

    \begin{minipage}{0.47\linewidth}
        \centering
        \includegraphics[clip=true,width=2.55in]{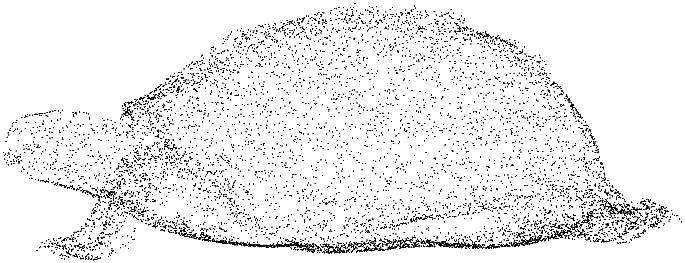}\\
        (a)
    \end{minipage}
    \hspace{0.025\linewidth}
    \begin{minipage}{0.47\linewidth}
        \centering
        \includegraphics[clip=true,width=2.55in]{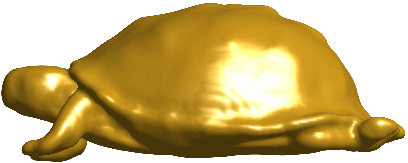}\\
        (b)
    \end{minipage}

    \vspace{0.5cm}

    \begin{minipage}{0.48\linewidth}
        \centering
        \includegraphics[clip=true,width=2.79in]{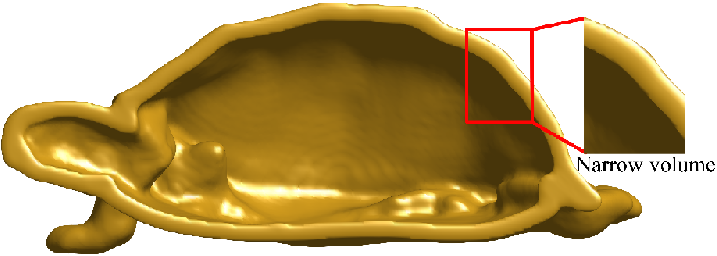}\\
        (c)
    \end{minipage}
    \hspace{0.025\linewidth}
    \begin{minipage}{0.47\linewidth}
        \centering
        \includegraphics[clip=true,width=2.38in]{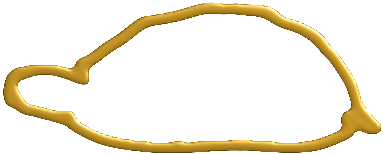}\\
        (d)
    \end{minipage}

    \caption{3D Reconstruction. The unorganized points, 3D object and narrow band are displayed.}
    \label{schfig1.eps}
\end{figure}

\begin{figure}[htbp]
\begin{minipage}{1\linewidth}
\centering
\includegraphics[clip=true,height=1.2in]{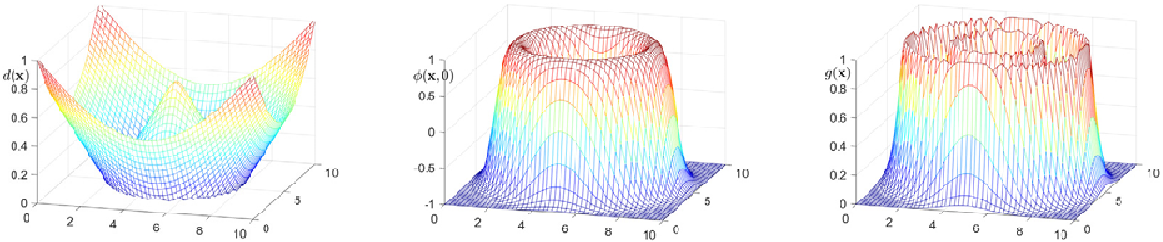}\\
\vspace{2mm}
\end{minipage}
\begin{minipage}{1\linewidth}
\centering
\includegraphics[clip=true,height=1.5in]{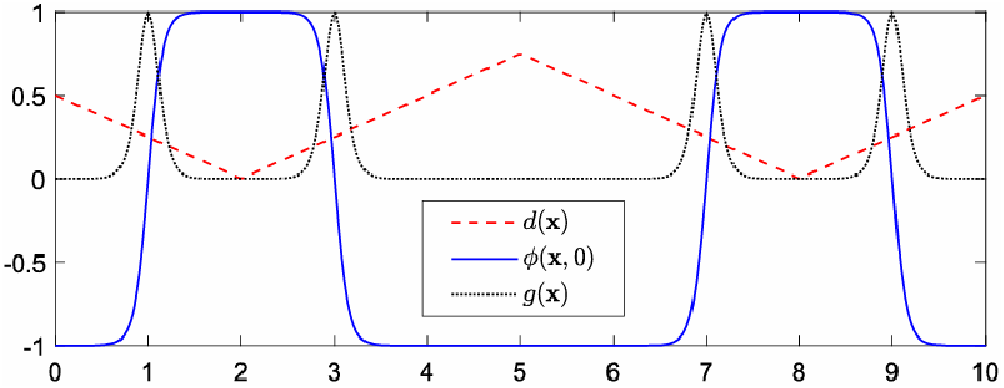}\\
\end{minipage}
    \caption{Schematic illustrations of signed distance function, phase-field function, and edge detection function. From the left to right in top row, the 2D illustrations of $d({\bf x})$, $\phi({\bf x},0)$, and $g({\bf x})$ are plotted. The bottom row shows the 1D cross sections. The figures were adapted from \cite{JJYangJCP} with the permission of Elsevier.} \label{fig2sch.eps}
\end{figure}

Within the spatial domain $\Omega$, the set of unorganized points is $S = \{{\bf X}_m \in \mathbb{R}^3 ~|~ m = 1, 2, ..., M_b\}$. We define the unsigned distance function as\ibe
d({\bf x}) = \min\limits_{1\le m \le M_b}|{\bf x}-{\bf X}_m|. \label{cddf1}\iee
Specifically, the distance function $d({\bf x})$ represents the unsigned distance from any point ${\bf x}$ in the domain $\Omega$ to its nearest neighbor ${\bf X}_m$ in the set $S$. An edge detection function (or control function) $g = g({\bf x})$ is introduced. We modify the AC model to the following form\ibe
\frac{\partial\phi}{\partial t} &=& -g\mu,\label{mac1}\\
\mu &=& \frac{F'(\phi)}{\epsilon^2} - \Delta\phi. \label{mac2}\iee
Here, ${\bf n}\cdot\nabla\phi|_{\partial\Omega} = 0$ is used on domain boundary. We have $\phi({\bf x},0) = \tanh\left(  \frac{\gamma - d({\bf x})}{\sqrt{2}\chi} \right)$. The value of $\chi$ is smaller than $\epsilon$. Then, we define $g({\bf x}) = 1-\phi^2({\bf x},0) + \lambda$. Here, $1-\phi^2({\bf x},0)$ acts as an interface indicator that concentrates the weighted gradient flow near the diffuse interface, while the added $\lambda$ guarantees strict positivity of $g$ (i.e., $g({\bf x})\ge\lambda$), preventing degeneracy and improving numerical robustness. We set $\lambda=10^{-10}$, which is sufficiently small so as not to affect the interface-localized dynamics. Figure \ref{schfig1.eps} provides schematic diagrams. In the top row of Fig. \ref{fig2sch.eps}, 2D schematic illustrations of distance function, phase-field function, and edge detection function are plotted. The bottom row of Fig. \ref{fig2sch.eps} plots the 1D cross sections of distance function, phase-field function, and edge detection function. From the definition of $g({\bf x})$, it is easy to observe the value of $g({\bf x})$ is bounded by $\lambda$ and $1+\lambda$.

We emphasize that $t$ is a pseudo-time for the $L^2$-gradient flow of $E(\phi)$ in Eq.~\eqref{Oene1}. The AC flow \eqref{ac1}--\eqref{ac2} and its edge-weighted form \eqref{mac1}--\eqref{mac2} satisfy the dissipation law $\tfrac{d}{dt}E(\phi) = -\|\sqrt{g}\mu\|^{2} \le 0$ (Theorem~2). This monotone dissipation provides a mathematically controlled relaxation toward steady states, in contrast to level‑set schemes in which periodic reinitialization of the signed‑distance field may perturb the interface and lacks a compatible discrete energy law \cite{OsherFedkiw2003,Li2010}. The reconstruction is taken at the steady state of this relaxation, where
\[ \partial_t \phi = 0 \;\Rightarrow\; g\,\mu = 0 \;\stackrel{g\ge \lambda>0}{\Rightarrow}\; \mu = 0 \;\Rightarrow\; \frac{F'(\phi)}{\epsilon^{2}} - \Delta \phi = 0, \]
under the homogeneous-Neumann boundary condition. The surface is extracted as the zero level set of the relaxed solution, $\{\mathbf{x}\in\Omega \mid \phi(\mathbf{x},t^{\ast})=0\}$. In practice, $g(\mathbf{x})$ localizes the evolution to a narrow band around the data-driven initial interface, and the unconditionally stable, decoupled Crank--Nicolson scheme in Section~\ref{sec3} ensures robust convergence with large time steps. Thus, introducing $\partial/\partial t$ provides a principled relaxation mechanism to compute the stationary contour.

{\bf Theorem 2}.{\bf 1}. The proposed model, under the homogeneous-Neumann boundary condition, maintains energy stability (i.e., the dissipation of the energy).\\
{\bf Proof}. By multiplying Eq. \eqref{mac1} with $\mu$ and implementing the integration, we derive\ibe
\left(\frac{\partial \phi}{\partial t},\mu\right) = -(g\mu,\mu) = -(\sqrt{g}\mu,\sqrt{g}\mu) = -\|\sqrt{g}\mu\|^2. \label{proff1}\iee
Multiplying Eq. \eqref{mac2} with $\partial\phi/\partial t$ and implementing the integration, we derive\ibe
\left(\mu,\frac{\partial\phi}{\partial t}\right) &=& \left( \frac{F'(\phi)}{\epsilon^2},\frac{\partial\phi}{\partial t}\right) - \left(\Delta\phi,\frac{\partial\phi}{\partial t}\right)\nonumber\\
&=& \frac{d}{dt}\int_\Omega \frac{F(\phi)}{\epsilon^2} ~d{\bm x} + \underbrace{\int_{\partial\Omega} \left( {\bf n}\cdot\nabla\phi \right)\frac{\partial\phi}{\partial t} ~ds}_{I} + \int_\Omega \nabla\phi \frac{\partial(\nabla\phi)}{\partial t} ~d{\bm x}\nonumber\\
&=& \frac{d}{dt}\int_\Omega \frac{F(\phi)}{\epsilon^2} ~d{\bm x}
+ \frac{d}{dt}\int_\Omega \frac{1}{2}|\nabla\phi|^2 ~d{\bm x}
= \frac{d}{dt}\int_\Omega \left[ \frac{F(\phi)}{\epsilon^2} + \frac{1}{2}|\nabla\phi|^2 \right] ~d{\bm x}. \label{proff2}\iee
As a result, term I can be eliminated since \( \mathbf{n} \cdot \nabla \phi |_{\partial \Omega} = 0 \). By integrating Eqs. \eqref{proff1} and \eqref{proff2}, we obtain\ibe
\frac{d}{dt}E(\phi) = \frac{d}{dt}\int_\Omega \left[ \frac{F(\phi)}{\epsilon^2} + \frac{1}{2}|\nabla\phi|^2 \right] ~d{\bm x}
= -\|\sqrt{g}\mu\|^2 \le 0. \iee
This derived inequality signifies that energy does not increase over time. Energy stability is a fundamental characteristic of the current phase-field model for volume reconstruction.

\section{Numerical scheme}\label{sec3}

We introduce a time-dependent variable $Q$ which theoretically equals to $1$. The equivalent governing equations read as
\begin{equation} \label{AC1}
    \frac{\partial \phi}{\partial t} = -g(x)\mu,
    \end{equation}
    \begin{equation} \label{AC2}
    \mu = -Q \frac{F'(\phi)}{\epsilon^2} + \Delta \phi,
    \end{equation}
    \begin{equation} \label{Qdf}
    \frac{d}{dt} \int_{\Omega} F(\phi) \, dx = Q \int_{\Omega}  \frac{\partial \phi}{\partial t} F'(\phi)\, dx.
    \end{equation}


\noindent\textit{Interpretation of $Q$.}
The weighted AC flow Eqs.~\eqref{mac1}--\eqref{mac2} is the \(L^2_g\)-steepest descent of \(E(\phi)\), where
\[
\langle u,v\rangle_g := \int_\Omega g(x)\,u\,v\,dx \;=\; (\sqrt{g}\,u,\sqrt{g}\,v).
\]
The time-dependent scalar \(Q=Q(t)\) introduced in Eqs.~\eqref{AC1}--\eqref{Qdf} is \emph{not} a new physical field; it is a Lagrange multiplier that enforces the continuous chain rule of the nonlinear energy under our semi-implicit discretization (cf.\ the discrete constraint \eqref{cnQdf} and the framework in [33]). Indeed, from \eqref{Qdf} and the chain rule we have
\[
\frac{d}{dt}\!\int_\Omega F(\phi)\,dx \;=\; \big(F'(\phi),\,\phi_t\big) \;=\; Q(t)\,\big(F'(\phi),\,\phi_t\big),
\]
so that
\[
Q(t)\;:=\;\frac{\frac{d}{dt}\int_\Omega F(\phi)\,dx}{\int_\Omega \phi_t\,F'(\phi)\,dx}\;\equiv\;1
\quad\text{whenever } \big(F'(\phi),\phi_t\big)\neq 0.
\]
At steady states \(\phi_t=0\), Eq. \eqref{AC1} implies \(g\,\mu=0\) and hence \(\mu=0\) (since \(g\ge\lambda>0\)), which reduces to the standard Euler–Lagrange equation of \(E\),
\[
\mu=0\ \Longleftrightarrow\ \frac{F'(\phi)}{\epsilon^2}-\Delta\phi=0.
\]
Consequently, \(Q\) does not change the equilibria or the zero level-set geometry of the reconstruction; it only rescales the \emph{temporal parametrization} of the descent in the metric \(\langle\cdot,\cdot\rangle_g\). In the fully discrete scheme, \(Q^{n+\frac12}\) is the unique scalar determined by \eqref{cnQdf} so that the original energy \(E(\phi)\) remains the Lyapunov functional while the update Eqs.~\eqref{cnAC1}--\eqref{cnAC2} stays linear and decoupled; thus \(Q\) serves as a structure-preserving device that encodes the exact discrete bookkeeping of the nonlinear energy, rather than a geometric regularizer.

\subsection{Discrete approach}
The 3D computational domain is $\Omega = (0,l_x)\times(0,l_y)\times(0,l_z)$. The uniform spatial step $h$ is defined as $h = l_x/N_x = l_y/N_y = l_z/N_z$. The numbers of grid along $x$-, $y$-, and $z$-directions are $N_x$, $N_y$, and $N_z$. The cubic cells are used to discretize $\Omega$. The center position of a cell is given by ${\bf x}_{ijk} = (x_i,y_j,z_k) = ((i-0.5)h,(j-0.5)h,(k-0.5)h)$. where $1 \leq i \leq N_x$, $1 \leq j \leq N_y$, and $1 \leq k \leq N_z$. Let $(\cdot)_{ijk}^n$ represent the numerical approximation of a variable at $(x_i, y_j, z_k, n\Delta t)$, where $\Delta t = \frac{T}{N_t}$ is the uniform time step, $T$ is the total computational time, and $N_t$ is the number of time iterations. The following discrete homogeneous-Neumann boundary conditions are applied to the discrete domain boundary $\partial\Omega_d$:\be
\phi_{0,jk} &=& \phi_{1,jk}, ~\phi_{N_x+1,jk}=\phi_{N_x,jk}, ~\mu_{0,jk}=\mu_{1,jk}, ~\mu_{N_x+1,jk} = \mu_{N_x,jk},\\
\phi_{i,0,k} &=& \phi_{i,1,k}, ~\phi_{i,N_y+1,k}=\phi_{i,N_y,k}, ~\mu_{i,0,k}=\mu_{i,1,k}, ~\mu_{i,N_y+1,k} = \mu_{i,N_y,k},\\
\phi_{ij,0} &=& \phi_{ij,1}, ~\phi_{ij,N_z+1} = \phi_{ij,N_z}, ~\mu_{ij,0}=\mu_{ij,1}, ~\mu_{ij,N_z+1} = \mu_{ij,N_z}.\ee
To evaluate the stability of energy in a discrete context, we define \be
(\phi,\psi)_h &=& h^3\sum\limits^{N_x}_{i=1}\sum\limits^{N_y}_{j=1}\sum\limits^{N_z}_{k=1} \phi_{ijk}\psi_{ijk},\\
\left(\nabla_d\phi,\nabla_d\psi\right)_e &=& h^3\left(\sum\limits^{N_x}_{i=0}\sum\limits^{N_y}_{j=1}\sum\limits^{N_z}_{k=1}
D_x\phi_{i+\frac{1}{2},jk}D_x\psi_{i+\frac{1}{2},jk} + \sum\limits^{N_x}_{i=1}\sum\limits^{N_y}_{j=0}\sum\limits^{N_z}_{k=1}
D_y\phi_{i,j+\frac{1}{2},k}D_y\psi_{i,j+\frac{1}{2},k} \right.\\
&&\left.+ \sum\limits^{N_x}_{i=1}\sum\limits^{N_y}_{j=1}\sum\limits^{N_z}_{k=0}
D_z\phi_{ij,k+\frac{1}{2}}D_z\psi_{ij,k+\frac{1}{2}} \right),\ee
We next define the discrete norm as
\begin{align*}
    \|\phi\|^2_h &= (\phi, \phi)_h, \quad \|\nabla_d \phi\|^2_e = (\nabla_d \phi, \nabla_d \phi)_e.
\end{align*}
The discrete derivatives are computed as
\begin{align*}
    D_x \phi_{i+\frac{1}{2},j,k} &= \frac{\phi_{i+1,j,k} - \phi_{i,j,k}}{h}, \quad
    D_y \phi_{i,j+\frac{1}{2},k} = \frac{\phi_{i,j+1,k} - \phi_{i,j,k}}{h}, \quad
    D_z \phi_{i,j,k+\frac{1}{2}} = \frac{\phi_{i,j,k+1} - \phi_{i,j,k}}{h}.
\end{align*}
We similarly define other quantities. The summation by parts is given by\be
&&\left(\phi,\Delta_d\psi\right)_h = h^3\sum\limits^{N_x}_{i=1}\sum\limits^{N_y}_{j=1}\sum\limits^{N_z}_{k=1} \phi_{ijk}\Delta_d\psi_{ijk}
= h^3\sum\limits^{N_x}_{i=1}\sum\limits^{N_y}_{j=1}\sum\limits^{N_z}_{k=1}\left(  \phi_{ijk}\frac{\psi_{i+1,jk}-\psi_{ijk}}{h^2} - \phi_{ijk}\frac{\psi_{ijk}-\psi_{i-1,jk}}{h^2}\right.\\
&&\left.+ \phi_{ijk}\frac{\psi_{i,j+1,k}-\psi_{ijk}}{h^2} - \phi_{ijk}\frac{\psi_{ijk}-\psi_{i,j-1,k}}{h^2}
+ \phi_{ijk}\frac{\psi_{ij,k+1}-\psi_{ijk}}{h^2} - \phi_{ijk}\frac{\psi_{ijk}-\psi_{ij,k-1}}{h^2}\right)\\
&=& h^3\left[ \sum\limits^{N_x}_{i=0}\sum\limits^{N_y}_{j=1}\sum\limits^{N_z}_{k=1}\left(\frac{\phi_{ijk}D_x\psi_{i+\frac{1}{2},jk}}{h}
-\frac{\phi_{i+1,jk}D_x\psi_{i+\frac{1}{2},jk}}{h}\right)
+ \sum\limits^{N_x}_{i=1}\sum\limits^{N_y}_{j=0}\sum\limits^{N_z}_{k=1}\left(\frac{\phi_{ijk}D_y\psi_{i,j+\frac{1}{2},k}}{h}
-\frac{\phi_{i,j+1,k}D_y\psi_{i,j+\frac{1}{2},k}}{h}\right) \right.\\
&&\left.+ \sum\limits^{N_x}_{i=1}\sum\limits^{N_y}_{j=1}\sum\limits^{N_z}_{k=0}\left(\frac{\phi_{ijk}D_z\psi_{ij,k+\frac{1}{2}}}{h}
-\frac{\phi_{ij,k+1}D_z\psi_{ij,k+\frac{1}{2}}}{h}\right)  \right]\\
&=& -h^3\left[ \sum\limits^{N_x}_{i=0}\sum\limits^{N_y}_{j=1}\sum\limits^{N_z}_{k=1}D_x\phi_{i+\frac{1}{2},jk}D_x\psi_{i+\frac{1}{2},jk}
+ \sum\limits^{N_x}_{i=1}\sum\limits^{N_y}_{j=0}\sum\limits^{N_z}_{k=1}D_y\phi_{i,j+\frac{1}{2},k}D_y\psi_{i,j+\frac{1}{2},k}
+ \sum\limits^{N_x}_{i=1}\sum\limits^{N_y}_{j=1}\sum\limits^{N_z}_{k=0}D_z\phi_{ij,k+\frac{1}{2}}D_z\psi_{ij,k+\frac{1}{2}}  \right]\\
&=& -\left(\nabla_d\phi,\nabla_d\psi\right)_e. \ee
To achieve second-order temporal accuracy, a CN scheme is adopted for discretizing the time derivatives, with the resulting scheme for Eqs. \eqref{AC1}, \eqref{AC2} and \eqref{Qdf} is written as
\begin{equation} \label{cnAC1}
    \frac{\phi^{n+1}_{ijk} - \phi^n_{ijk}}{\Delta t} = -g_{ijk}\mu^{n+\frac{1}{2}}_{ijk},
\end{equation}
\begin{equation} \label{cnAC2}
    \mu^{n+\frac{1}{2}}_{ijk} = Q^{n+\frac{1}{2}} \frac{F'(\phi^*_{ijk})}{\epsilon^2} - \Delta_d\left( \frac{\phi^{n+1}_{ijk} + \phi^n_{ijk}}{2} \right) + S \left( \frac{\phi^{n+1}_{ijk} + \phi^n_{ijk}}{2} - \phi^*_{ijk} \right),
\end{equation}
\begin{equation} \label{cnQdf}
    \left(F(\phi^{n+1}) - F(\phi^n), {\bf 1}\right)_h = Q^{n+\frac{1}{2}} \left(F'(\phi^*),\phi^{n+1} - \phi^n\right)_h.
\end{equation}
where $\Delta_d\phi_{ijk} = (\phi_{i+1,jk}+\phi_{i-1,jk}+\phi_{i,j+1,k}+\phi_{i,j-1,k}+\phi_{ij,k+1}+\phi_{ij,k-1}-6\phi_{ijk})/h^2$,  $(\cdot)^* = \frac{3}{2}(\cdot)^n - \frac{1}{2}(\cdot)^{n-1}$, i.e., $\phi^* := \tfrac{3}{2}\phi^n - \tfrac{1}{2}\phi^{n-1}$, $S > 0$ is a stabilizing constant. We note that the coupling between \( \phi^{n+1} \) and \( Q^{n+\frac{1}{2}} \) makes the solving process inefficient.
To decouple the unknown Lagrange multiplier $Q^{n+\frac{1}{2}}$ from the update of $\phi$, we split the new solution and the chemical potential at each grid point $(i,j,k)$ into a $Q$-independent part $(\cdot)_1$ and a part linear in $Q$ $(\cdot)_2$. Specifically, we define $\phi_{ijk}^{n+1}=\phi_{1,ijk}^{n+1}+Q^{n+\frac{1}{2}}\phi_{2,ijk}^{n+1}$ and $\mu_{ijk}^{n+\frac{1}{2}}=\mu_{1,ijk}^{n+\frac{1}{2}}+Q^{n+\frac{1}{2}}\mu_{2,ijk}^{n+\frac{1}{2}}$. Here, $\phi_{1,ijk}$ and $\mu_{1,ijk}$ collect the contributions that do not depend on $Q$ (the homogeneous subproblem), while $\phi_{2,ijk}$ and $\mu_{2,ijk}$ represent the response to the nonlinear source term $F'(\phi^*)/\epsilon^2$ under the same linear operator. This affine decomposition allows us to solve two linear subproblems first and then determine $Q^{n+\frac{1}{2}}$ from the discrete energy constraint \eqref{cnQdf}, thereby removing the coupling and improving efficiency.
By substituting these expressions into \eqref{cnAC1} and \eqref{cnAC2}, we have \ibe
\frac{\phi_{1,ijk}^{n+1} + Q^{n+\frac{1}{2}} \phi_{2,ijk}^{n+1} - \phi^n_{ijk}}{\Delta t} &=& -g_{ijk}\left( \mu_{1,ijk}^{n+\frac{1}{2}} + Q^{n+\frac{1}{2}} \mu_{2,ijk}^{n+\frac{1}{2}} \right),\label{lCNs1}\\
\mu_{1,ijk}^{n+\frac{1}{2}} + Q^{n+\frac{1}{2}} \mu_{2,ijk}^{n+\frac{1}{2}} &=& Q^{n+\frac{1}{2}} \frac{F'(\phi_{ijk}^*)}{\epsilon^2} - \frac{\Delta_d}{2} \left( \phi_{1,ijk}^{n+1} + Q^{n+\frac{1}{2}} \phi_{2,ijk}^{n+1} + \phi_{ijk}^n \right) \nonumber\\
&&+ S \left( \frac{1}{2} \left( \phi_{1,ijk}^{n+1} + Q^{n+\frac{1}{2}} \phi_{2,ijk}^{n+1} + \phi^n_{ijk} \right) - \phi^*_{ijk} \right). \label{lCNs2}\iee
By using a splitting strategy, we decompose Eqs. \eqref{lCNs1} and \eqref{lCNs2} into
\begin{align}
    \frac{\phi_{1,ijk}^{n+1} - \phi_{ijk}^n}{\Delta t} &= -g_{ijk}\mu_{1,ijk}^{n+ \frac{1}{2}} \label{eq:update_phi1}, \\
    \mu_{1,ijk}^{n+ \frac{1}{2}} &= - \frac{\Delta_d}{2} \left( \phi_{1,ijk}^{n+1} + \phi_{ijk}^n \right)
    + S \left( \frac{1}{2} \left( \phi_{1,ijk}^{n+1} + \phi^n_{ijk} \right) - \phi^*_{ijk} \right) \label{eq:update_mu1},
\end{align}
and
\begin{align}
    \frac{\phi_{2,ijk}^{n+1}}{\Delta t} &= -g_{ijk}\mu_{2,ijk}^{n+ \frac{1}{2}} \label{eq:update_phi2}, \\
    \mu_{2,ijk}^{n+ \frac{1}{2}} &= \frac{F'(\phi^*_{ijk})}{\epsilon^2} - \frac{\Delta_d}{2} \phi_{2,ijk}^{n+1} + \frac{S}{2} \phi_{2,ijk}^{n+1} \label{eq:update_mu2}.
\end{align}
We use Eqs. \eqref{eq:update_phi1} and \eqref{eq:update_mu1} to update $\phi_{1,ijk}^{n+1}$, and Eqs. \eqref{eq:update_phi2} and \eqref{eq:update_mu2} to update $\phi_{2,ijk}^{n+1}$. Substituting $\phi^{n+1}_{ijk} = \phi_{1,ijk}^{n+1} + Q^{n + \frac{1}{2}} \phi_{2,ijk}^{n+1}$ into Equation \eqref{cnQdf}, we obtain
\noindent
\begin{equation}
\begin{aligned}
\left( \left[ \frac{1}{4} \left( \left( \phi_1^{n+1} + Q^{n+\frac{1}{2}} \phi_2^{n+1} \right)^2 - 1 \right)^2 - F(\phi^n_{ijk}) \right],{\bf 1}\right)_h
= Q^{n+\frac{1}{2}}\left(F'(\phi^*), \left( \phi_{1,ijk}^{n+1} + Q^{n+\frac{1}{2}} \phi_{2,ijk}^{n+\frac{1}{2}} - \phi^n_{ijk} \right)\right)_h
.\end{aligned}
\label{lCNs3}
\end{equation}
Here, ${\bf 1}$ represents the vector with all entries equal to $1$. Substitute the computed $\phi_{1,ijk}^{n+1}$ and $\phi_{2,ijk}^{n+1}$ into equation \eqref{lCNs3}, leaving  $Q^{n + \frac{1}{2}}$ as the only unknown variable in the equation. By rearranging equation \eqref{lCNs2}, we obtain
\begin{align*}
    \quad & \left(\frac{1}{4} \left( \left( \phi_1^{n+1} \right)^2 + \left( Q^{n+\frac{1}{2}} \phi_2^{n+1} \right)^2 + 2 Q^{n+\frac{1}{2}} \phi_1^{n+1} \phi_2^{n+1} - 1 \right)^2 - F(\phi^n),{\bf 1} \right)_h \notag \\
    &= Q^{n+\frac{1}{2}} \left(\phi_1^{n+1},F'(\phi^*)\right)_h + \left(Q^{n+\frac{1}{2}}\right)^2\left(F'(\phi^*),\phi_2^{n+\frac{1}{2}}\right)_h - Q^{n+\frac{1}{2}}\left(F'(\phi^*),\phi^n \right)_h.
\end{align*}

    \begin{align*}
        \Rightarrow \quad & \left(\frac{1}{4} \left( \phi_1^{n+1} \right)^4 + \frac{1}{4}\left( \phi_2^{n+1} \right)^4 \left( Q^{n+\frac{1}{2}} \right)^4  + \left( \phi_1^{n+1} \right)^3 \phi_2^{n+1} Q^{n+\frac{1}{2}} + \phi_1^{n+1} \left( \phi_2^{n+1} \right)^3 \left( Q^{n+\frac{1}{2}} \right)^3\right.\nonumber \\
        & \left.\quad + \frac{3}{2} \left( \phi_1^{n+1} \right)^2 \left( \phi_2^{n+1} \right)^2 \left( Q^{n+\frac{1}{2}} \right)^2 - \frac{1}{2} \left( \phi_1^{n+1} \right)^2 - \frac{1}{2} \left( \phi_2^{n+1} \right)^2 \left( Q^{n+\frac{1}{2}} \right)^2\right.\nonumber \\
        & \left.\quad - \phi_1^{n+1} \phi_2^{n+1} Q^{n+\frac{1}{2}} + \frac{1}{4} - F(\phi^n),{\bf 1}\right)_h \nonumber \\
        &= \left( Q^{n+\frac{1}{2}} \right)^2 \left(\phi_2^{n+1},F'(\phi^*)\right)_h + Q^{n+\frac{1}{2}} \left(\phi_1^{n+1},F'(\phi^*)\right)_h - Q^{n+\frac{1}{2}}\left(F'(\phi^*),\phi^n\right)_h.
    \end{align*}
By expanding the terms and simplifying, we arrive at
\begin{align}
    & \left( Q^{n+\frac{1}{2}} \right)^4 \left(\frac{1}{4} \left( \phi_2^{n+1} \right)^4,{\bf 1}\right)_h + \left( Q^{n+\frac{1}{2}} \right)^3 \left(\phi_1^{n+1},\left( \phi_2^{n+1} \right)^3 \right)_h \nonumber \\
    & \quad + \left( Q^{n+\frac{1}{2}} \right)^2 \left(\frac{3}{2} \left( \phi_1^{n+1} \right)^2,\left( \phi_2^{n+1} \right)^2\right)_h
    - \left( Q^{n+\frac{1}{2}} \right)^2\left(\frac{1}{2} \left( \phi_2^{n+1} \right)^2,{\bf 1}\right)_h -
    \left( Q^{n+\frac{1}{2}} \right)^2\left(\phi_2^{n+1},F'(\phi^*) \right)_h\nonumber \\
    & \quad + Q^{n+\frac{1}{2}} \left(\left( \phi_1^{n+1} \right)^3,\phi_2^{n+1}\right)_h - Q^{n+\frac{1}{2}}\left(\phi_1^{n+1},\phi_2^{n+1}\right)_h - Q^{n+\frac{1}{2}}\left(F'(\phi^*),\phi^{n+1}_1\right)_h + Q^{n+\frac{1}{2}}\left(F'(\phi^*),\phi^n \right)_h \nonumber \\
    & \quad + \left(\frac{1}{4} \left( \phi_1^{n+1} \right)^4 - \frac{1}{2} \left( \phi_1^{n+1} \right)^2 - F(\phi^n) + \frac{1}{4},{\bf 1}\right)_h = 0.
    \label{simplified}
\end{align}
By applying the Newton iteration method to Eq. \eqref{simplified}, we iteratively solve for $Q^{n+\frac{1}{2}}$  until the change in $Q^{n+\frac{1}{2}}$ is less than the specified tolerance, thereby obtaining the solution for $Q^{n+\frac{1}{2}}$. Once $Q^{n+\frac{1}{2}}$, $\phi_{1,ijk}^{n+1}$ and $\phi_{2,ijk}^{n+1}$ are computed, we can then update $\phi^{n+1}_{ijk}$ using $\phi^{n+1}_{ijk} = \phi_{1,ijk}^{n+1} + Q^{n + \frac{1}{2}} \phi_{2,ijk}^{n+1}$.

\subsection{Estimation of discrete stability}\label{subsec:discrete_stability}
In this subsection, we implement estimation of discrete energy stability. Multiplying Eq. \eqref{cnAC1} with $\Delta t\mu^{n+\frac{1}{2}}_{ijk}$ and implementing the integration, we get\ibe
\left(\phi^{n+1}-\phi^n,\mu^{n+\frac{1}{2}}\right)_h = -\|\sqrt{g}\mu^{n+\frac{1}{2}}\|^2_h. \label{dproff3}\iee
Multiplying Eq. \eqref{cnAC2} with $\phi^{n+1}_{ijk}-\phi^n_{ijk}$ and implementing the integration, we get\ibe
\left(\mu^{n+\frac{1}{2}},\phi^{n+1}-\phi^n\right)_h &=& Q^{n+\frac{1}{2}}\left(\frac{F'(\phi^*)}{\epsilon^2},\phi^{n+1}-\phi^n\right)_h
- \left(\Delta_d\left( \frac{\phi^{n+1}+\phi^n}{2}\right),\phi^{n+1}-\phi^n\right)_h \nonumber\\
&&+ S\left(\frac{\phi^{n+1}+\phi^n}{2}-\phi^*,\phi^{n+1}-\phi^n\right)_h, \label{dproff4}\iee
where \be
-\left(\Delta_d\left( \frac{\phi^{n+1}+\phi^n}{2}\right),\phi^{n+1}-\phi^n\right)_h = \frac{1}{2}\left(\nabla_d\phi^{n+1}-\nabla_d\phi^n,\nabla_d\phi^{n+1}+\nabla_d\phi^n\right)_e
 = \frac{1}{2}\left( \|\nabla_d\phi^{n+1}\|^2_e - \|\nabla_d\phi^n\|^2_e\right), \ee
and\be
S\left(\frac{\phi^{n+1}+\phi^n}{2}-\phi^*,\phi^{n+1}-\phi^n\right)_h &=&
S\left(\frac{\phi^{n+1}+\phi^n}{2} - \left( \frac{3}{2}\phi^n - \frac{1}{2}\phi^{n-1} \right),\phi^{n+1}-\phi^n\right)_h \nonumber\\
&=& \frac{S}{2}\left(\phi^{n+1}+\phi^n,\phi^{n+1}-\phi^n\right)_h - \frac{S}{2}\left(3\phi^n-\phi^{n-1},\phi^{n+1}-\phi^n\right)_h\nonumber\\
&=& \frac{S}{2}\left(\|\phi^{n+1}\|^2_h - \|\phi^n\|^2_h\right) - \frac{S}{2}(\|\phi^{n+1}\|^2_h - \|\phi^n\|^2_h \nonumber\\
&=& - \frac{1}{2}(\|\phi^{n+1}-\phi^n\|^2_h - \|\phi^n-\phi^{n-1}\|^2_h + \|\phi^{n+1}-2\phi^n+\phi^{n-1}\|^2_h))\nonumber\\
&=& \frac{S}{4}\left( \|\phi^{n+1}-\phi^n\|^2_h - \|\phi^n-\phi^{n-1}\|^2_h + \|\phi^{n+1}-2\phi^n+\phi^{n-1}\|^2_h  \right). \ee
Therefore, we have\ibe
\left(\mu^{n+\frac{1}{2}},\phi^{n+1}-\phi^n\right)_h &=& Q^{n+\frac{1}{2}}\left(\frac{F'(\phi^*)}{\epsilon^2},\phi^{n+1}-\phi^n\right)_h
+ \frac{1}{2}\left( \|\nabla_d\phi^{n+1}\|^2_e - \|\nabla_d\phi^n\|^2_e\right)\nonumber\\
&&+ \frac{S}{4}\left( \|\phi^{n+1}-\phi^n\|^2_h - \|\phi^n-\phi^{n-1}\|^2_h + \|\phi^{n+1}-2\phi^n+\phi^{n-1}\|^2_h  \right). \label{proff3}\iee
Multiplying $\frac{1}{\epsilon^2}$ with Eq. \eqref{cnQdf}, we have\ibe
\frac{1}{\epsilon^2}\left(F(\phi^{n+1}) - F(\phi^n), {\bf 1}\right)_h = Q^{n+\frac{1}{2}} \left(\frac{F'(\phi^*)}{\epsilon^2},\phi^{n+1} - \phi^n\right)_h.\label{proff4}\iee
Combining the results in Eqs. \eqref{proff1}, \eqref{proff3}, and \eqref{proff4}, we derive\ibe
&&\frac{1}{\epsilon^2}\left(F(\phi^{n+1}) - F(\phi^n), {\bf 1}\right)_h + \frac{1}{2}\left( \|\nabla_d\phi^{n+1}\|^2_e - \|\nabla_d\phi^n\|^2_e\right)
+ \frac{S}{4}\left( \|\phi^{n+1}-\phi^n\|^2_h - \|\phi^n-\phi^{n-1}\|^2_h\right)\nonumber\\
&=& -\|\sqrt{g}\mu^{n+\frac{1}{2}}\|^2_h - \frac{S}{4}\|\phi^{n+1}-2\phi^n+\phi^{n-1}\|^2_h \le 0. \label{eq:discrete_energy_dissipation} \iee
For clarity and completeness, we define the discrete energy functional, which is consistent with the stability estimate in Eq.~\eqref{eq:discrete_energy_dissipation}:
\ibe\label{eq:discrete_energy}
\tilde{E}(\phi^{n+1},\phi^n) := \frac{1}{\epsilon^2}\left(F(\phi^{n+1}),{\bf 1}\right)_h + \frac{1}{2}\|\nabla_d\phi^{n+1}\|^2_e + \frac{S}{4}\|\phi^{n+1}-\phi^n\|^2_h.\iee
In particular, the discrete dissipation relation implies that $\tilde{E}(\phi^{n+1},\phi^n) \le \tilde{E}(\phi^{n},\phi^{n-1})$.
The proof of discrete energy dissipation law is completed.

\subsection{Numerical implementation}
This subsection provides a concise overview of the numerical implementation. To simplify the expressions, we begin by reformulating Eqs. \eqref{eq:update_phi1} and \eqref{eq:update_mu1} as\ibe
L\left(\phi^{n+1}_{1,ijk}~,~\mu^{n+\frac{1}{2}}_{1,ijk}\right) = \left(d^n_{1,ijk}~,~d^n_{2,ijk}\right), \label{0CN1}\iee
where the operator \( L \) and the right-hand side terms are defined as  \be
L\left(\phi^{n+1}_{1,ijk}~,~\mu^{n+\frac{1}{2}}_{1,ijk}\right) = \left( \frac{\phi^{n+1}_{1,ijk}}{\Delta t} + g_{ijk}\mu^{n+\frac{1}{2}}_{1,ijk}~,
~\mu^{n+\frac{1}{2}}_{1,ijk}-\frac{S}{2}\phi^{n+1}_{1,ijk}+\frac{1}{2}\Delta_d\phi^{n+1}_{1,ijk}\right),\ee
and \be
d^n_{1,ijk} &=& \frac{\phi^n_{1,ijk}}{\Delta t}, \\
d^n_{2,ijk} &=& S\left( - \phi^n_{1,ijk} + \frac{1}{2} \phi^{n-1}_{1,ijk} \right) - \frac{1}{2}\Delta_d\phi^n_{ijk}.\ee

Similarly, we reformulate Eqs. \eqref{eq:update_phi2} and \eqref{eq:update_mu2} as\ibe
L\left(\phi^{n+1}_{2,ijk}~,~\mu^{n+\frac{1}{2}}_{2,ijk}\right) = \left(p^n_{1,ijk}~,~p^n_{2,ijk}\right), \label{0CN2}\iee
with \be
L\left(\phi^{n+1}_{2,ijk}~,~\mu^{n+\frac{1}{2}}_{2,ijk}\right) = \left( \frac{\phi^{n+1}_{2,ijk}}{\Delta t} + g_{ijk}\mu^{n+\frac{1}{2}}_{2,ijk}~,
~\mu^{n+\frac{1}{2}}_{2,ijk}-\frac{S}{2}\phi^{n+1}_{2,ijk}+\frac{1}{2}\Delta_d\phi^{n+1}_{2,ijk}\right),\ee
and \be
p^n_{1,ijk} &=& 0, \\
p^n_{2,ijk} &=& \frac{F'(\phi^*_{ijk})}{\epsilon^2}.\ee

The Gauss--Seidel (GS) iteration is used to solve the following $2 \times 2$ matrix equations until specified tolerances are met. Specifically, we first solve the matrix equation for $\phi^{n+1,m+1}_{1}$ until its convergence criterion is satisfied, and then solve the matrix equation for $\phi^{n+1,m+1}_{2}$ until its convergence criterion is met. i.e.,

\begin{equation}
\begin{gathered}
\begin{pmatrix}
\frac{1}{\Delta t} & g_{ijk} \\
-\frac{S}{2}-\frac{3}{h^2} & 1
\end{pmatrix}
\begin{pmatrix}
\phi^{n+1,m+1}_{1,ijk} \\
\mu^{n+\frac{1}{2}}_{1,ijk}
\end{pmatrix} \\
= \begin{pmatrix}
d^n_{1,ijk} \\
d^n_{2,ijk} - \frac{1}{2h^2}\left( \phi^{n+1,m}_{1,i+1,jk}+\phi^{n+1,m}_{1,i-1,jk}+\phi^{n+1,m}_{1,i,j+1,k}+\phi^{n+1,m}_{1,i,j-1,k}
+ \phi^{n+1,m}_{1,ij,k+1}+\phi^{n+1,m}_{1,ij,k-1} \right)
\end{pmatrix},
\end{gathered}
\end{equation}

\begin{equation}
\begin{gathered}
\begin{pmatrix}
\frac{1}{\Delta t} & g_{ijk} \\
-\frac{S}{2}-\frac{3}{h^2} & 1
\end{pmatrix}
\begin{pmatrix}
\phi^{n+1,m+1}_{2,ijk} \\
\mu^{n+\frac{1}{2}}_{2,ijk}
\end{pmatrix} \\
= \begin{pmatrix}
p^n_{1,ijk} \\
p^n_{2,ijk} - \frac{1}{2h^2}\left( \phi^{n+1,m}_{2,i+1,jk}+\phi^{n+1,m}_{2,i-1,jk}+\phi^{n+1,m}_{2,i,j+1,k}+\phi^{n+1,m}_{2,i,j-1,k}
+ \phi^{n+1,m}_{2,ij,k+1}+\phi^{n+1,m}_{2,ij,k-1} \right)
\end{pmatrix}.
\end{gathered}
\end{equation}

We use $m+1$ and $m$ to represent the solutions after and before the $(m+1)$-th iteration. The stop criteria are defined as $\|\phi^{n+1,m+1}_1 - \phi^{n+1,m}_1\|^2_h \le 1$e-$6$ for $\phi_1$, and $\|\phi^{n+1,m+1}_2 - \phi^{n+1,m}_2\|^2_h \le 1$e-$6$ for $\phi_2$. Once these criteria are satisfied, we update the solutions by setting $\phi^{n+1}_{1,ijk} = \phi^{n+1,m+1}_{1,ijk}$ and $\phi^{n+1}_{2,ijk} = \phi^{n+1,m+1}_{2,ijk}$. Finally, we compute $Q^{n+\frac{1}{2}}$ using Newton iteration, thereby obtaining $\phi^{n+1}_{ijk}$ based on our proposed scheme. Please refer to Algorithm 1 for the implementation.
\begin{algorithm}
\KwData{Calculated variables at previous time levels}
Calculate $\phi^{n+1}_{1,ijk}$ and $\mu^{n+\frac{1}{2}}_{1,ijk}$ by solving Eqs. \eqref{eq:update_phi1} and \eqref{eq:update_mu1} with a GS-type iterative algorithm\;
\While{$\|\phi^{n+1,m+1}_{1}-\phi^{n+1,m}_{1}\|^2_h > 1$\mbox{e}-$6$}
{$(\phi^{n+1,m+1}_{1,ijk},\mu^{n+1,m+1}_{1,ijk}) = \mbox{GS-cycle}(\phi^{n+1,m}_{1,ijk},d_1^n,d_2^n,u)$\;
$\star$ $u$ denotes the total number of GS relaxations\;
}
Update $\phi^{n+1}_{2,ijk}$ and $\mu^{n+\frac{1}{2}}_{2,ijk}$ by solving Eqs. \eqref{eq:update_phi2} and \eqref{eq:update_mu2} with a GS-type iterative algorithm\;
\While{$\|\phi^{n+1,m+1}_{2}-\phi^{n+1,m}_{2}\|^2_h > 1$\mbox{e}-$6$}
{$(\phi^{n+1,m+1}_{2,ijk},\mu^{n+1,m+1}_{2,ijk}) = \mbox{GS-cycle}(\phi^{n+1,m}_{2,ijk},q_1^n,q_2^n,v)$\;
$\star$ $v$ denotes the total number of GS relaxations\;
}
Compute $Q^{n+\frac{1}{2}}$ by solving Eq. \eqref{simplified} using Newton iteration algorithm with a tolerance of $1$e-$6$\;
Update $\phi^{n+1}_{ijk}$\ using $\phi^{n+1}_{ijk} = \phi_{1,ijk}^{n+1} + Q^{n + \frac{1}{2}} \phi_{2,ijk}^{n+1}$\;
Compute the discrete energy $\tilde{E}(\phi^{n+1},\phi^n)$ defined in Eq.~\eqref{eq:discrete_energy}\;
\KwResult{Output the updated phase-field function $\phi^{n+1}_{ijk}$}
\caption{Implementation of computation}
\end{algorithm}

\vspace{1\baselineskip}
{\bf Remark 3}.{\bf 1}. We adopt a uniform Cartesian grid with the standard 7-point finite-difference stencil for three main reasons. First, on Cartesian meshes the discrete gradient–divergence pair satisfies an exact summation-by-parts identity, which directly underpins the discrete energy-dissipation law in Subsection \ref{subsec:discrete_stability} and \eqref{eq:discrete_energy_dissipation}–\eqref{eq:discrete_energy}. Second, the scheme is robust and reproducible: it avoids mesh-quality sensitivities, is straightforward to implement, and parallelizes efficiently due to its strictly local stencil. Third, in our formulation the edge-weighted flow $\partial_t\phi=-g\,\mu$ concentrates the dynamics within a narrow interfacial band, so although the grid is uniform, the effective computational work is largely localized near the surface, which alleviates the perceived inefficiency of a uniform mesh. These considerations motivate the use of uniform-grid FDM as a clean, structure-preserving baseline for our reconstruction experiments.

\section{Numerical simulations}\label{sec4}
In this section, extensive numerical results are shown to verify the accuracy, stability, and capability of our method. In all simulations, we set $S = \frac{2}{\epsilon^2}$; the reason for this choice is discussed in Subsection~\ref{subsec:S_effect}.

\subsection{Stability of energy law}
\begin{figure}[htbp]
    \centering

    \begin{minipage}{0.32\linewidth}
        \centering
        \includegraphics[clip=true,width=1.5in]{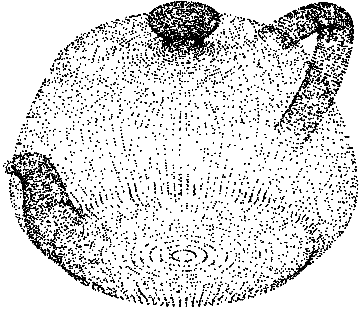}
        \\(a)\vspace{2mm}
    \end{minipage}
    \begin{minipage}{0.32\linewidth}
        \centering
        \includegraphics[clip=true,width=1.5in]{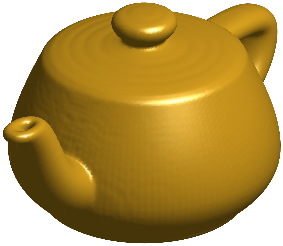}
        \\(b)\vspace{2mm}
    \end{minipage}
    \begin{minipage}{0.32\linewidth}
        \centering
        \includegraphics[clip=true,width=1.5in]{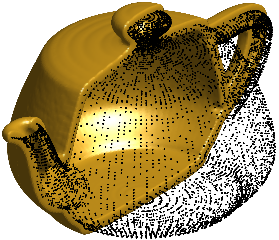}
        \\(c)\vspace{2mm}
    \end{minipage}
    \\[1mm]

    \begin{minipage}{0.32\linewidth}
        \centering
        \includegraphics[clip=true,width=1.7in]{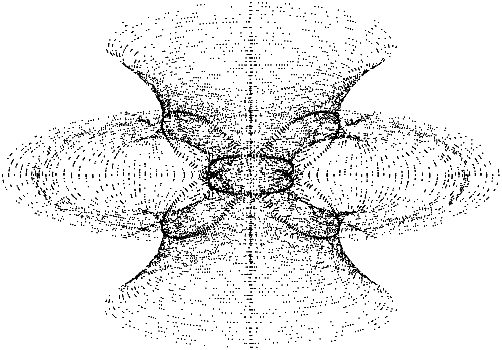}
        \\(d)\vspace{2mm}
    \end{minipage}
    \begin{minipage}{0.32\linewidth}
        \centering
        \includegraphics[clip=true,width=1.7in]{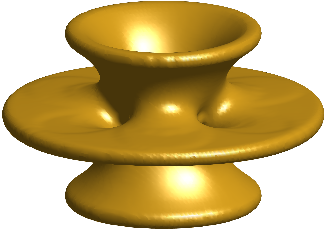}
        \\(e)\vspace{2mm}
    \end{minipage}
    \begin{minipage}{0.32\linewidth}
        \centering
        \includegraphics[clip=true,width=1.7in]{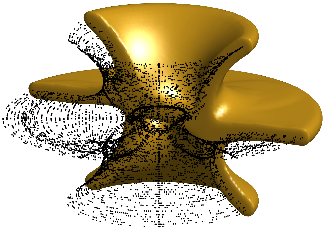}
        \\(f)\vspace{2mm}
    \end{minipage}
    \\[1mm]

    \begin{minipage}{0.32\linewidth}
        \centering
        \includegraphics[clip=true,width=1.5in]{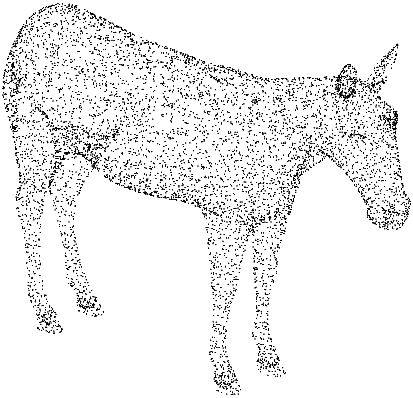}
        \\(g)\vspace{2mm}
    \end{minipage}
    \begin{minipage}{0.32\linewidth}
        \centering
        \includegraphics[clip=true,width=1.5in]{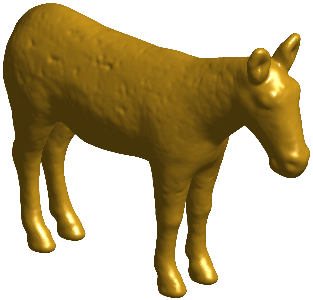}
        \\(h)\vspace{2mm}
    \end{minipage}
    \begin{minipage}{0.32\linewidth}
        \centering
        \includegraphics[clip=true,width=1.5in]{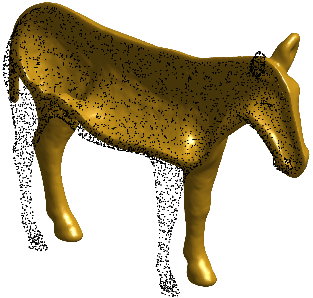}
        \\(i)\vspace{2mm}
    \end{minipage}

    \caption{3D narrow volume reconstructions of three different objects. }
    \label{enefig1}
\end{figure}

\begin{figure}[htbp]
    \centering

    \begin{minipage}{0.62\linewidth}
        \centering
        \includegraphics[width=\linewidth]{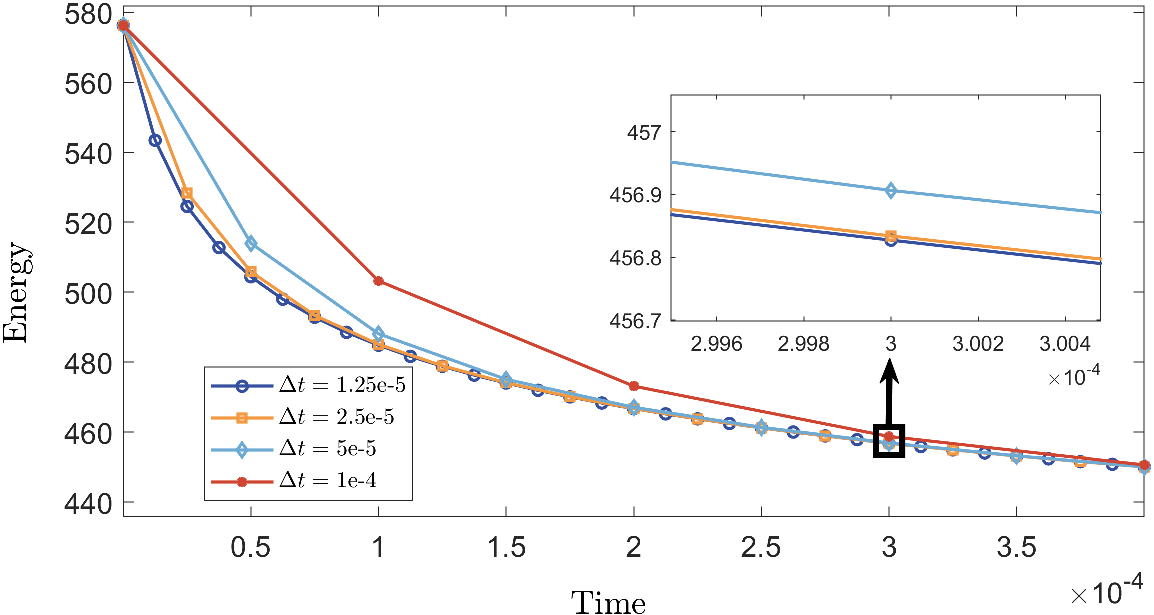}
        \\(a) Teapot
    \end{minipage}\\[2mm]
    \begin{minipage}{0.62\linewidth}
        \centering
        \includegraphics[width=\linewidth]{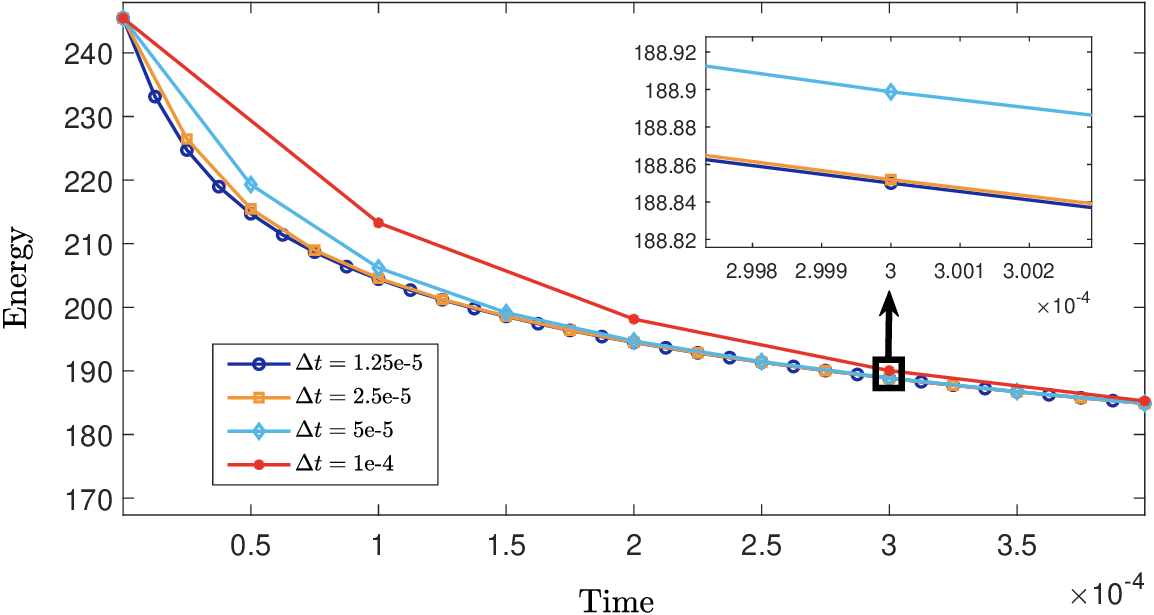}
        \\(b) Costa--Hoffman--Meeks
    \end{minipage}\\[2mm]
    \begin{minipage}{0.62\linewidth}
        \centering
        \includegraphics[width=\linewidth]{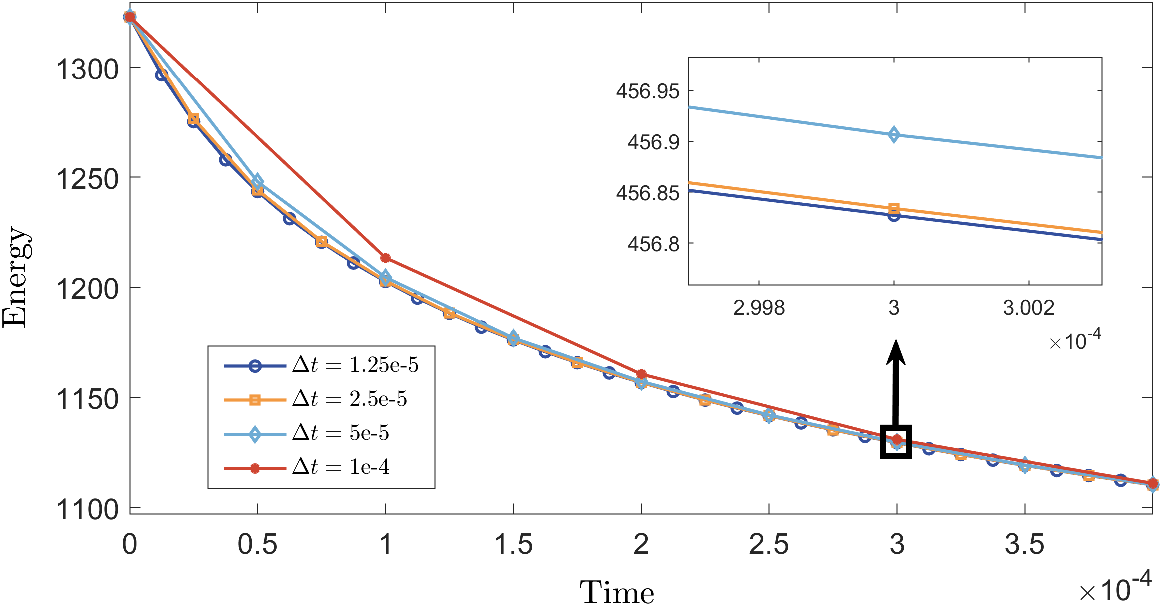}
        \\(c) Horse
    \end{minipage}\\[2mm]
    \caption{Energy stability for three different objects with different time steps.}
    \label{enefig2}
\end{figure}

We first preprocess the scattered point data using MATLAB and store the refined position data of the point cloud. The computational domain for the teapot is defined as $(0,2)\times(0,1)\times(0,1)$. We set $h = \frac{1}{128}$ and $\epsilon = 0.0113$. For the reconstruction of horse and Costa--Hoffman--Meeks surface, similar preprocessing steps are conducted using MATLAB. We use $h = 1/64$ and $\epsilon = 0.0450$. In all cases, simulations are carried out up to $t = 2$e-$4$. In Fig. \ref{enefig1}(a), (d) and (g), we display the points cloud. The reconstructed results with $\Delta t = 2\delta t$ for teapot and costa and with $\Delta t = 10\delta t$ for horse are shown. The energy curves with different time steps are displayed in Fig. \ref{enefig2}(a), (b) and (c). We observe that the energy dissipation law (i.e., energy stability) is satisfied.

\subsection{Accuracy experiment}
\begin{figure}[htbp]
    \centering
    \begin{minipage}{0.6\linewidth}
        \centering
        \includegraphics[clip=true,width=\linewidth]{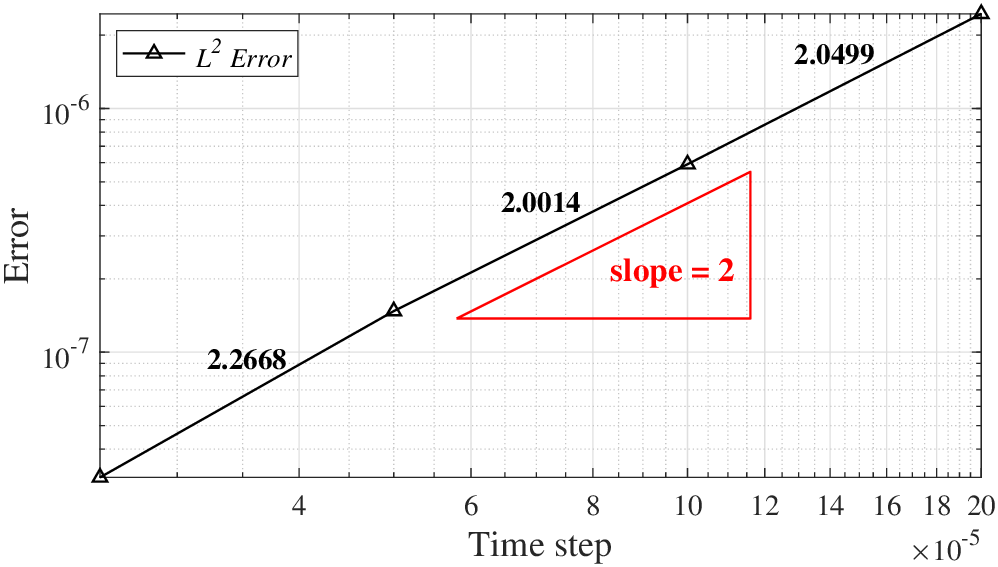}\\
        (a) Teapot
    \end{minipage}\\[2mm]
    \begin{minipage}{0.6\linewidth}
        \centering
        \includegraphics[clip=true,width=\linewidth]{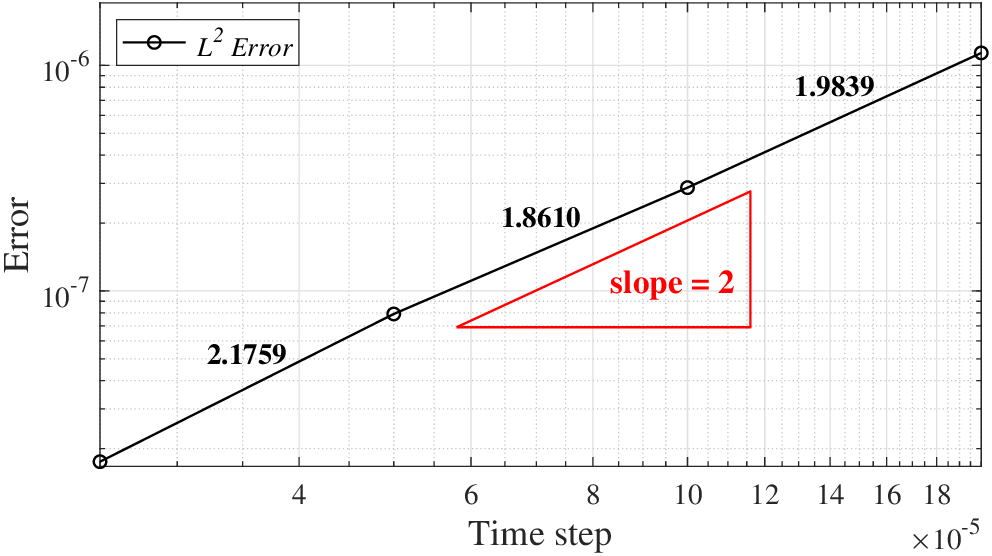}\\
        (b) Costa--Hoffman--Meeks
    \end{minipage}\\[2mm]
    \begin{minipage}{0.6\linewidth}
        \centering
        \includegraphics[clip=true,width=\linewidth]{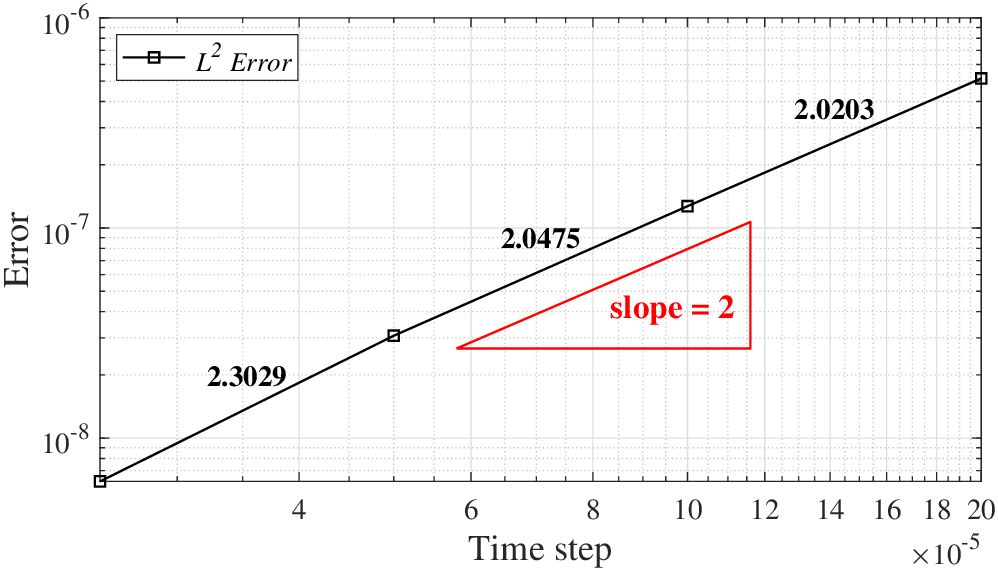}\\
        (c) Horse
    \end{minipage}
    \caption{Accuracy experiment of our method. In each figure, the errors versus time steps are plotted.}
    \label{acc1}
\end{figure}

Three cases used in last subsection are used to implement the accuracy test, all using $\epsilon = 0.0113$ and their own grid spacing $h$ as in the previous section. We obtain the reference solutions using a fine time step of $\hat{\delta t} = 1.25 \times 10^{-5}$. The $L^2$-errors are computed using increasingly coarser time steps: $\Delta t = 2\hat{\delta t}$, $4\hat{\delta t}$, $8\hat{\delta t}$, and $16\hat{\delta t}$. The results are plotted in Fig.~\ref{acc1}. The convergence rates between adjacent time steps are listed in the figures and are approximately $2$, indicating that our method achieves second-order accuracy in time. This is further supported by the observation that the slopes of the error curves are approximately equal to $2$. In addition, the $L^2$-errors and the convergence rates for the three datasets (teapot, Costa--Hoffman--Meeks, and horse) are also summarized in three tables, see Tables~\ref{tab:acc-teapot}--\ref{tab:acc-horse}.

\noindent\textit{Comparison with other methods.} We compare our scheme with other methods. Cai et al.\ propose a BDF2 time–marching strategy for the AC–based narrow-volume model\cite{Cai2025}; Kong et al.\  recast the nonlinear term via a time-dependent auxiliary variable and develop two lower-boundedness–preserving linear schemes (AV-A, AV-B) that are unconditionally energy-stable \cite{Kong2025}; All comparisons are on the Costa case. In Table~\ref{tab:method-comparison}, rows 1–2 list the $L^2$ errors and convergence rates of AV-A and AV-B, sourced from ~\cite{Kong2025}; row 3 lists the BDF2 results, sourced from ~\cite{Cai2025}; and the last row reports our Lagrange–multiplier (LM) scheme. For clarity, the baseline numbers are taken directly from the cited papers without recomputation. The comparison shows that the AV family is first-order accurate, whereas both BDF2 and our method achieve second-order accuracy with compatible rates. This cross-paper consistency supports the observed orders and motivates the use of second-order time discretizations when higher temporal accuracy is needed.

\begin{table}[htbp]
    \centering
    \caption{Teapot: $L^2$-errors and convergence rates with respect to different time steps.}
    \label{tab:acc-teapot}
    \begin{tabular}{llllllll}
    \toprule
    $\Delta t$ & $16\hat{\delta t}$ &  & $8\hat{\delta t}$ &  & $4\hat{\delta t}$ &  & $2\hat{\delta t}$ \\
    \midrule
    $L^2$-error & 1.135734e-6 &  & 2.871104e-7 &  & 7.903755e-8 &  & 1.749187e-8 \\
    Rate &  & 1.983947 &  & 1.860995 &  & 2.175854 &  \\
    \bottomrule
    \end{tabular}
    \end{table}

    \begin{table}[htbp]
    \centering
    \caption{Costa--Hoffman--Meeks: $L^2$-errors and convergence rates with respect to different time steps.}
    \label{tab:acc-costa}
    \begin{tabular}{llllllll}
    \toprule
    $\Delta t$ & $16\hat{\delta t}$ &  & $8\hat{\delta t}$ &  & $4\hat{\delta t}$ &  & $2\hat{\delta t}$ \\
    \midrule
    $L^2$-error & 2.444562e-6 &  & 5.903799e-7 &  & 1.474521e-7 &  & 3.063900e-8 \\
    Rate &  & 2.049860 &  & 2.001397 &  & 2.266805 &  \\
    \bottomrule
    \end{tabular}
    \end{table}

    \begin{table}[htbp]
    \centering
    \caption{Horse: $L^2$-errors and convergence rates with respect to different time steps.}
    \label{tab:acc-horse}
    \begin{tabular}{llllllll}
    \toprule
    $\Delta t$ & $16\hat{\delta t}$ &  & $8\hat{\delta t}$ &  & $4\hat{\delta t}$ &  & $2\hat{\delta t}$ \\
    \midrule
    $L^2$-error & 5.150218e-7 &  & 1.269541e-7 &  & 3.071002e-8 &  & 6.223618e-9 \\
    Rate &  & 2.020327 &  & 2.047526 &  & 2.302884 &  \\
    \bottomrule
    \end{tabular}
    \end{table}

    \begin{table}[htbp]
        \centering
        \caption{Costa--Hoffman--Meeks: $L^2$-errors and convergence rates (multi-method) with respect to different time steps.}
        \label{tab:method-comparison}
        \setlength{\tabcolsep}{4pt} 
        \small                      
        \resizebox{\linewidth}{!}{%
        \begin{tabular}{lllllllll}
          \toprule
          Method & $\Delta t$ & $16\hat{\delta t}$ & & $8\hat{\delta t}$ & & $4\hat{\delta t}$ & & $2\hat{\delta t}$ \\
          \midrule
          \multirow{2}{*}{AV-A} & $L^2$-error & 3.561456e-3 & & 1.680863e-3 & & 7.248435e-4 & & 2.423552e-4 \\
                                    & Rate         &                & 1.083265      &                & 1.213461      &               & 1.580547      & \\
          \addlinespace
          \multirow{2}{*}{AV-B} & $L^2$-error & 3.718156e-3 & & 1.800460e-3 & & 8.019209e-4 & & 2.779629e-4 \\
                                    & Rate         &                & 1.046222      &                & 1.166833      &               & 1.528568      & \\
          \addlinespace
          \multirow{2}{*}{BDF2} & $L^2$-error & 2.805997e-6 & & 6.450041e-7 & & 1.778174e-7 & & 4.150774e-8 \\
                                    & Rate         &               & 2.121118      &                & 1.858922      &                & 2.098949      & \\
          \addlinespace
          \multirow{2}{*}{\textbf{LM (Ours)}} & $L^2$-error & 1.135734e-6 & & 2.871104e-7 & & 7.903755e-8 & & 1.749187e-8 \\
                                    & Rate         &               & \textbf{1.983947}      &                & \textbf{1.860995}      &                & \textbf{2.175854}      & \\
          \bottomrule
        \end{tabular}%
        }
      \end{table}

\subsection{Point cloud-based narrow volume reconstruction}\label{vobjects}
\begin{figure}[htbp]
\begin{minipage}{0.33\linewidth}
\centering
\includegraphics[clip=true,width=1.2in]{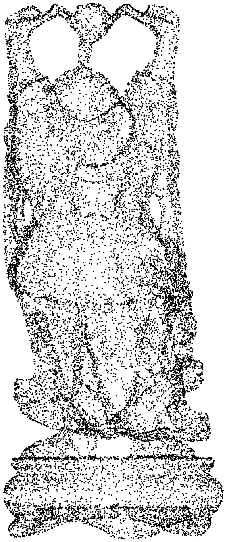}\\
\vspace{2mm}
\end{minipage}
\begin{minipage}{0.33\linewidth}
\centering
\includegraphics[clip=true,width=1.2in]{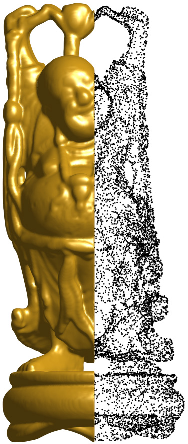}\\
\vspace{2mm}
\end{minipage}
\begin{minipage}{0.32\linewidth}
\centering
\includegraphics[clip=true,width=1.2in]{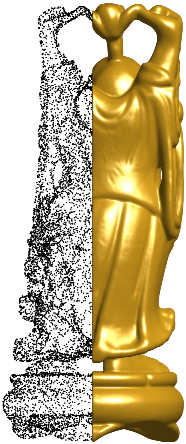}\\
\vspace{2mm}
\end{minipage}\\
\begin{minipage}{0.33\linewidth}
\centering
\includegraphics[clip=true,width=1.8in]{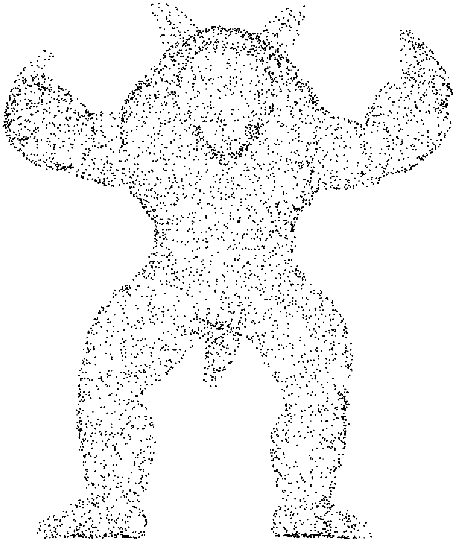}\\
\vspace{2mm}
\end{minipage}
\begin{minipage}{0.33\linewidth}
\centering
\includegraphics[clip=true,width=1.8in]{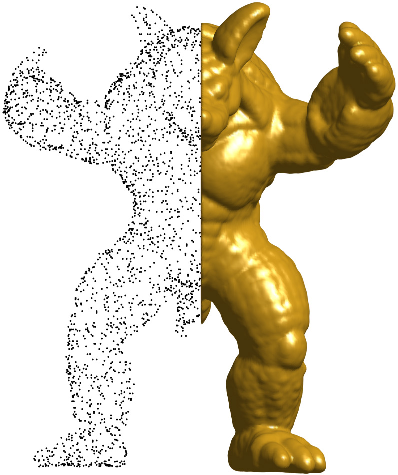}\\
\vspace{2mm}
\end{minipage}
\begin{minipage}{0.32\linewidth}
\centering
\includegraphics[clip=true,width=1.8in]{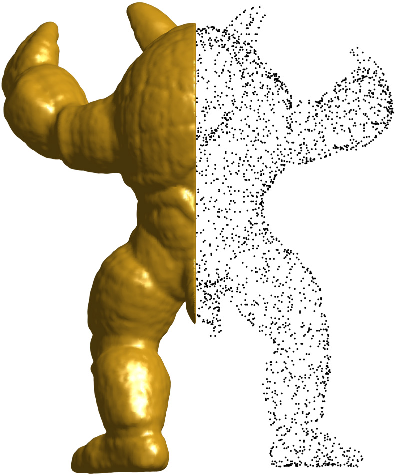}\\
\vspace{2mm}
\end{minipage}\\
\begin{minipage}{0.33\linewidth}
\centering
\includegraphics[clip=true,width=1.8in]{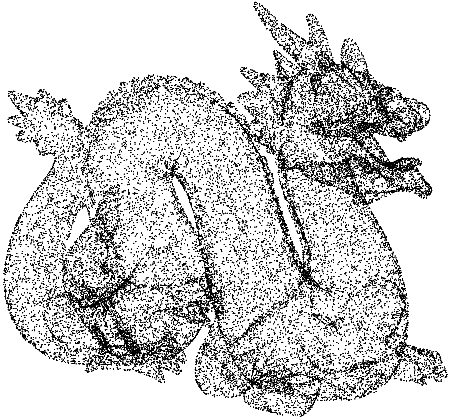}\\
\vspace{4mm}
\end{minipage}
\begin{minipage}{0.33\linewidth}
\centering
\includegraphics[clip=true,width=1.8in]{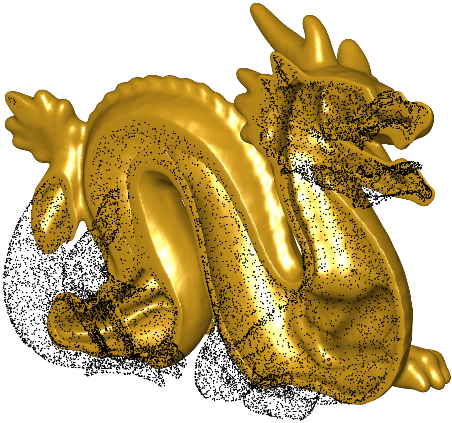}\\
\vspace{4mm}
\end{minipage}
\begin{minipage}{0.32\linewidth}
\centering
\includegraphics[clip=true,width=2in]{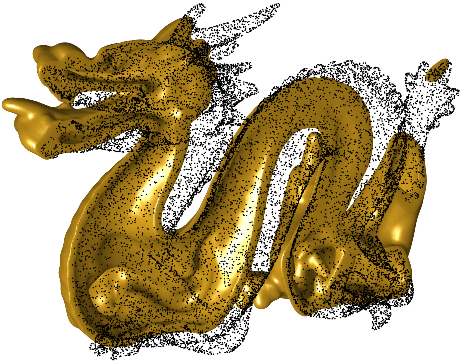}\\
\vspace{0mm}
\end{minipage}
\caption{Reconstructions of happy Buddha, Armadillo, and Stanford dragon from points cloud. } \label{multic1.eps}
\end{figure}

\begin{figure}[htbp]
\begin{minipage}{0.33\linewidth}
\centering
\includegraphics[clip=true,width=1.1in]{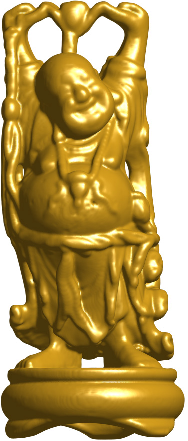}\\
\vspace{2mm}
\end{minipage}
\begin{minipage}{0.33\linewidth}
\centering
\includegraphics[clip=true,width=1.7in]{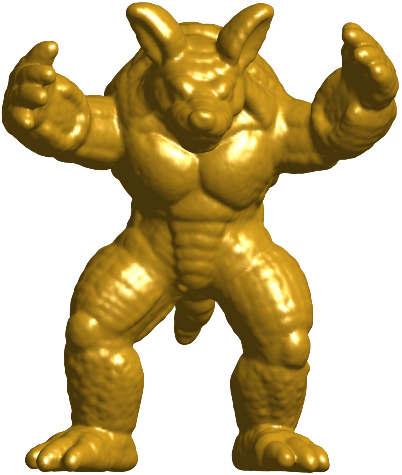}\\
\vspace{2mm}
\end{minipage}
\begin{minipage}{0.32\linewidth}
\centering
\includegraphics[clip=true,width=1.9in]{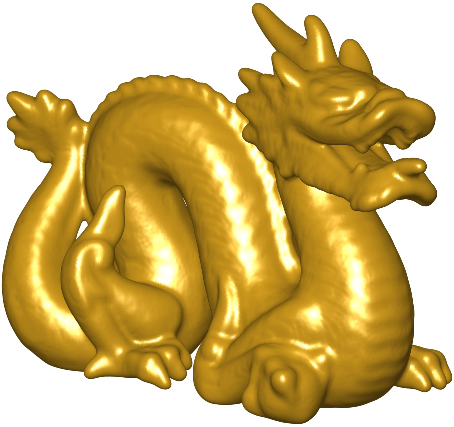}\\
\vspace{2mm}
\end{minipage}\\
\begin{minipage}{0.33\linewidth}
\centering
\includegraphics[clip=true,width=2in]{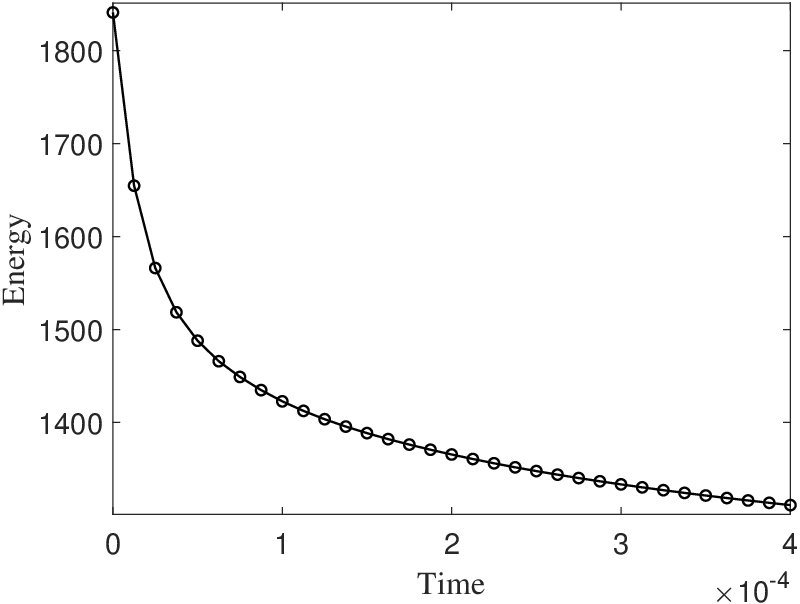}\\
\vspace{2mm}
\end{minipage}
\begin{minipage}{0.33\linewidth}
\centering
\includegraphics[clip=true,width=2in]{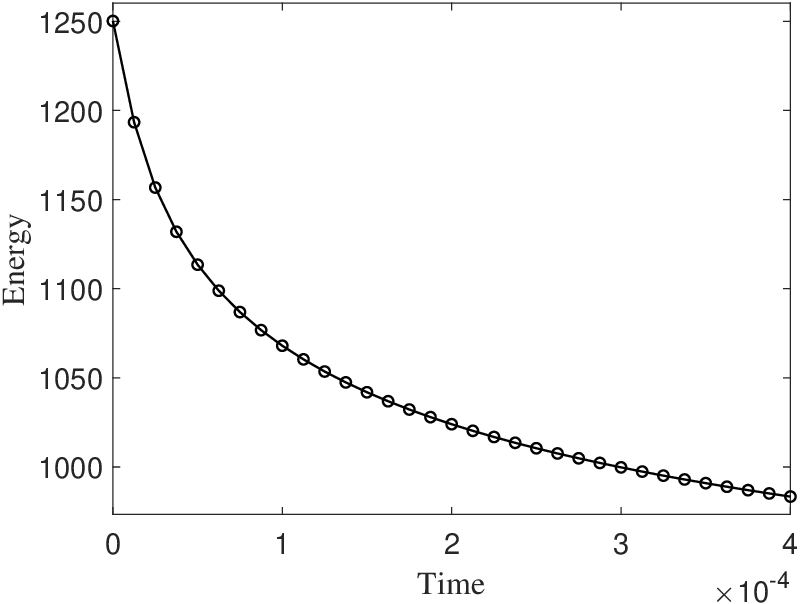}\\
\vspace{2mm}
\end{minipage}
\begin{minipage}{0.32\linewidth}
\centering
\includegraphics[clip=true,width=2in]{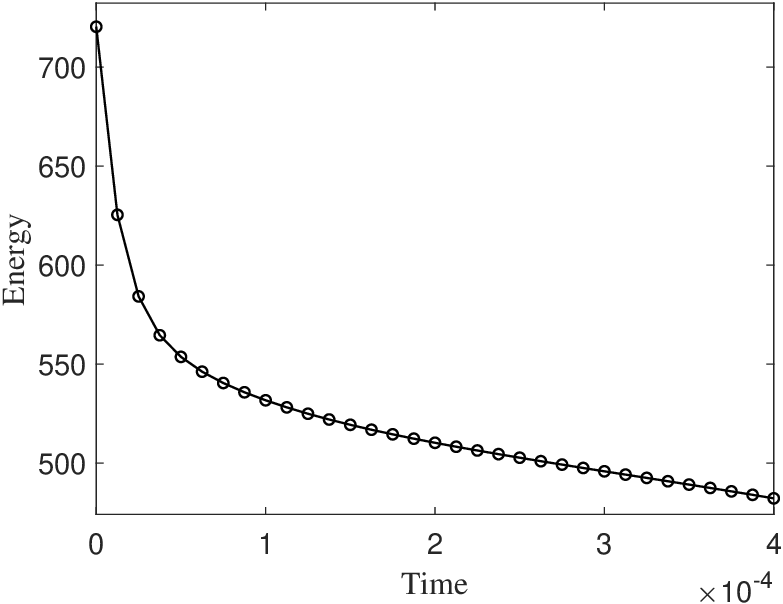}\\
\vspace{2mm}
\end{minipage}
\caption{Energy curves with respect to three different 3D volumes.} \label{multic2.eps}
\end{figure}

In this subsection, we aim to use our proposed algorithm to reconstruct various 3D objects. The domains are $(0,1.2)\times(0,2.4)\times(0,1.2)$ for a happy Buddha, $ (0,2.2)\times(0,2.5)\times(0,2)$ for an Armadillo, and $(0,1)\times(0,0.8)\times(0,0.6)$ for a Stanford dragon. We assign $\epsilon = 0.0068, ~0.0144, ~0.0036$, and $h = 0.3/64, ~1/100, ~1/400$ for different cases. The point sets and numerical results are displayed in Fig. \ref{multic1.eps}. Figure \ref{multic2.eps} presents three-dimensional profiles. The respective energy evolution curves over time are plotted.

\begin{figure}[htbp]
\begin{minipage}{0.33\linewidth}
\centering
\includegraphics[clip=true,width=2in]{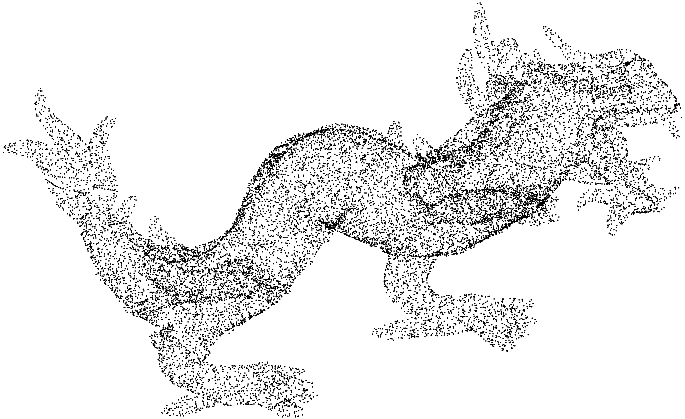}\\
\vspace{3mm}
\end{minipage}
\begin{minipage}{0.33\linewidth}
\centering
\includegraphics[clip=true,width=2in]{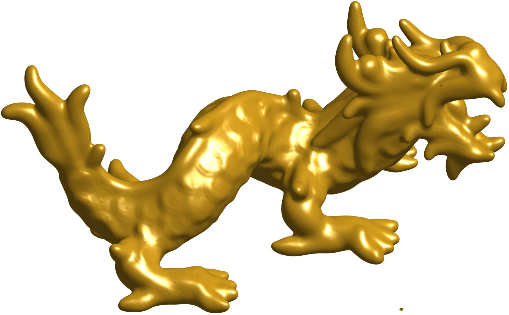}\\
\vspace{3mm}
\end{minipage}
\begin{minipage}{0.32\linewidth}
\centering
\includegraphics[clip=true,width=2in]{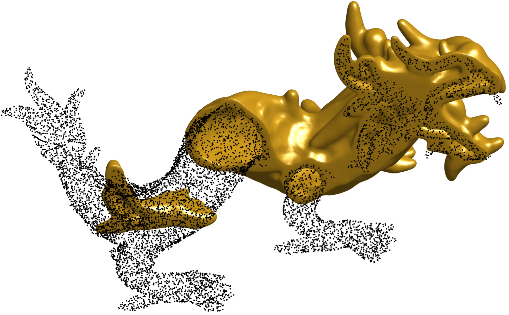}\\
\vspace{3mm}
\end{minipage}\\
\begin{minipage}{0.33\linewidth}
\centering
\includegraphics[clip=true,width=1.0in]{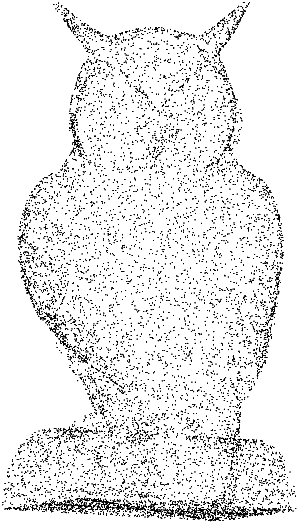}\\
\vspace{3mm}
\end{minipage}
\begin{minipage}{0.33\linewidth}
\centering
\includegraphics[clip=true,width=1.0in]{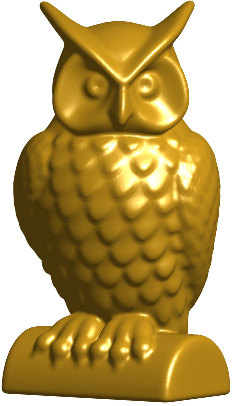}\\
\vspace{3mm}
\end{minipage}
\begin{minipage}{0.32\linewidth}
\centering
\includegraphics[clip=true,width=1.0in]{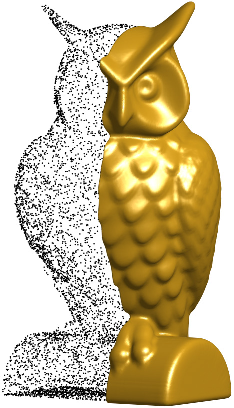}\\
\vspace{3mm}
\end{minipage}\\
\begin{minipage}{0.33\linewidth}
\centering
\includegraphics[clip=true,width=1.8in]{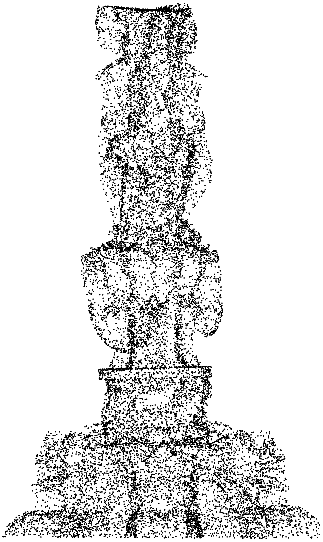}\\
\vspace{3mm}
\end{minipage}
\begin{minipage}{0.33\linewidth}
\centering
\includegraphics[clip=true,width=1.8in]{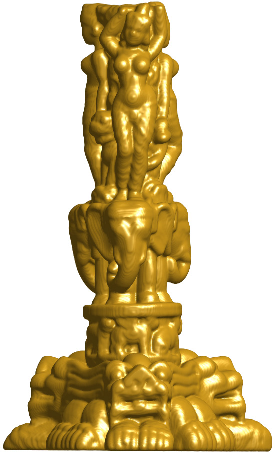}\\
\vspace{3mm}
\end{minipage}
\begin{minipage}{0.32\linewidth}
\centering
\includegraphics[clip=true,width=1.8in]{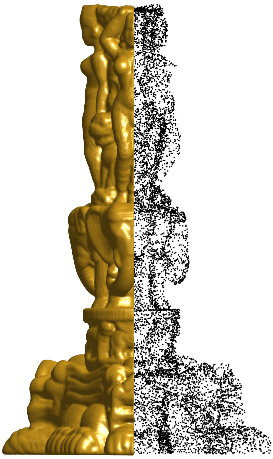}\\
\vspace{3mm}
\end{minipage}
\caption{Reconstructions of three different objects from points cloud.} \label{multic3.eps}
\end{figure}

Next, we implement simulations using points cloud of three different objects. To enhance the visualization of the scattered point sets, we show them in Fig. \ref{multic3.eps}, where only a portion of the points is shown for clarity. The reconstructed results are also shown in Fig. \ref{multic3.eps}. These visualizations demonstrate that the proposed algorithm effectively reconstructs a variety of 3D volumes with high fidelity.

\subsection{Effects of different transitional interface thickness parameters on reconstructions}
\begin{figure}[htbp]
    \centering

    \begin{minipage}{0.24\linewidth}
        \centering
        \includegraphics[width=0.98\linewidth]{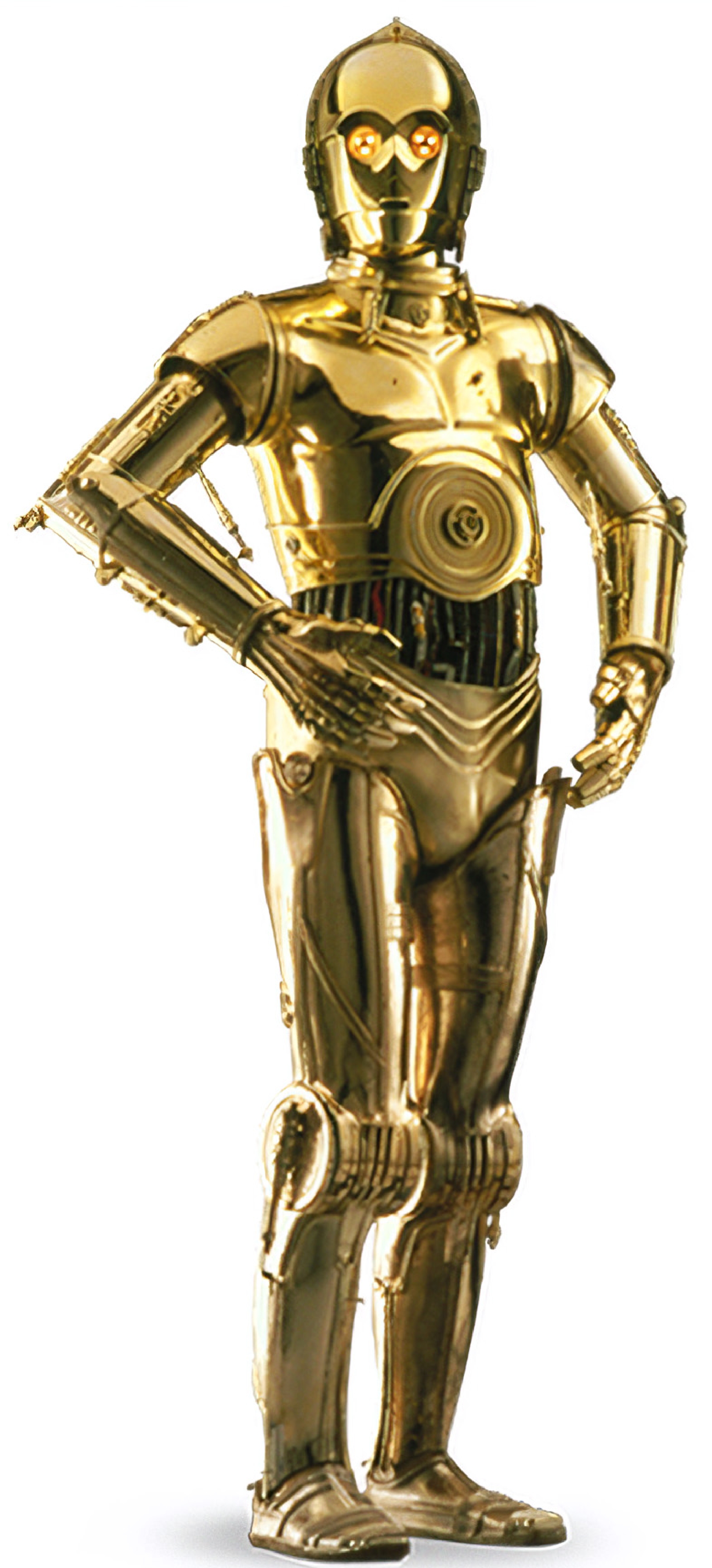}
    \end{minipage}
    \hfill
    \begin{minipage}{0.24\linewidth}
        \centering
        \includegraphics[width=\linewidth]{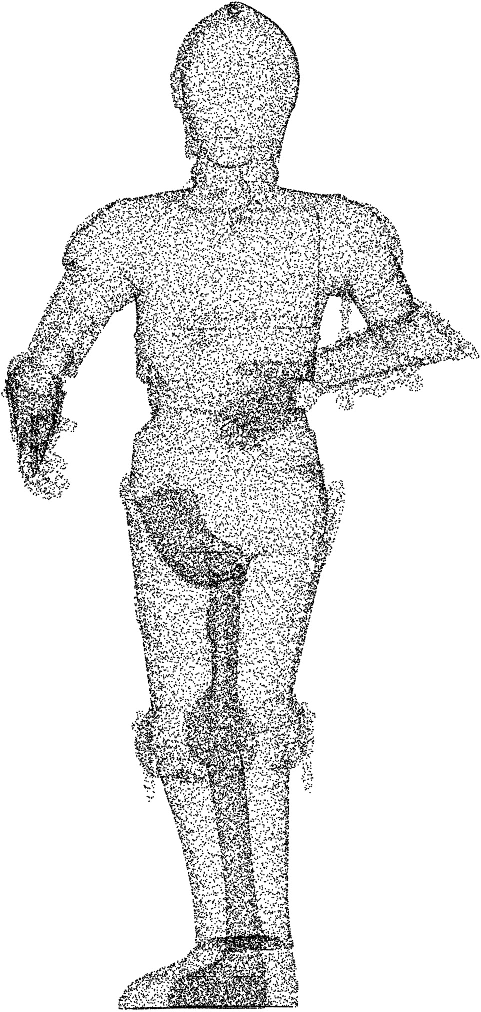}
    \end{minipage}
    \hfill
    \begin{minipage}{0.24\linewidth}
        \centering
        \includegraphics[width=\linewidth]{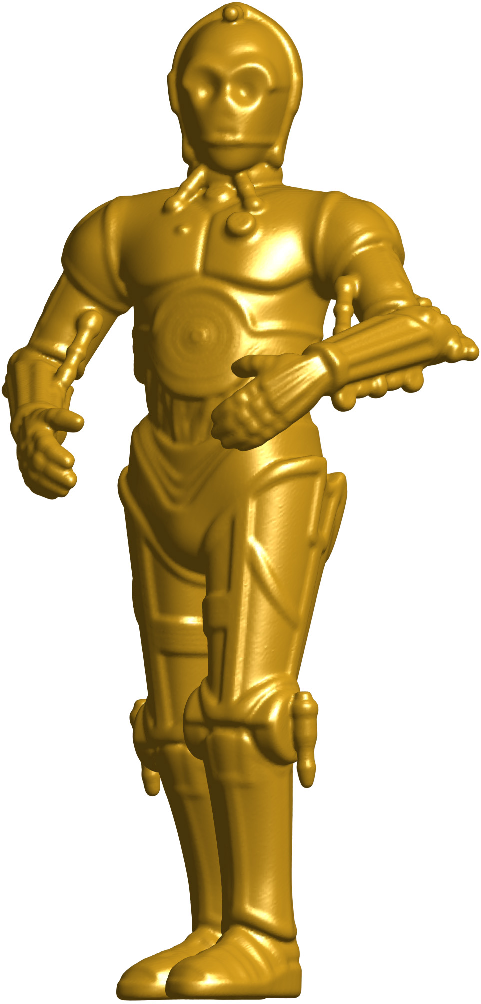}
    \end{minipage}
    \hfill
    \begin{minipage}{0.24\linewidth}
        \centering
        \includegraphics[width=\linewidth]{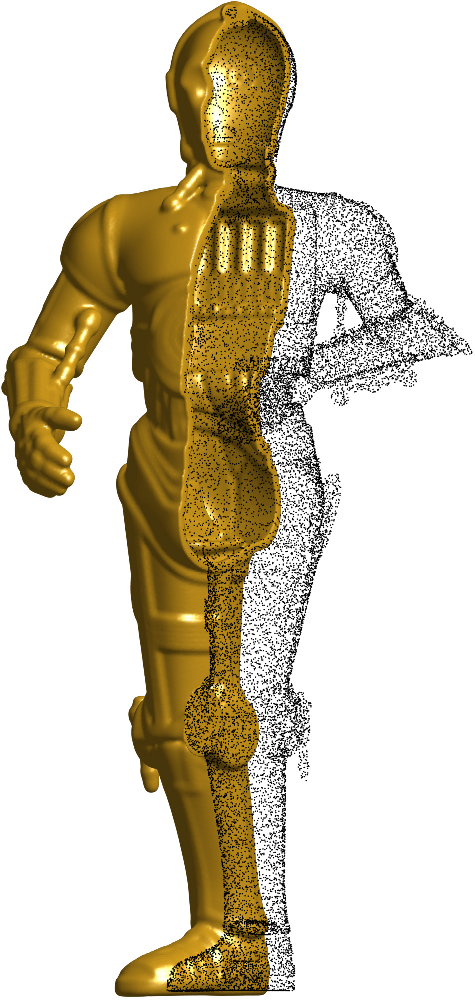}
    \end{minipage}

    \vspace{3mm}
    \caption{3D narrow volume reconstruction of C-3PO. From left to right: the appearance of C-3PO in the movie, the point cloud, the reconstructed volume, and a cut-off view of the reconstructed volume.} \label{c3po_volume}
\end{figure}

As stated in Section \ref{sec2}, the thickness of the transition layer of the reconstructed model is associated with a positive parameter $\epsilon$. Here, we set $\epsilon_a = 4$e-$4$, $\epsilon_b = 8$e-$4$, $\epsilon_c = 1.6$e-$3$, and $\epsilon_d = 3.2$e-$3$. In the test sequence, each $\epsilon$ is approximately twice the preceding one. We use C-3PO as an experimental sample, who is a robot from the renowned movie series \textit{Star Wars}. After preprocessing the 3D model file obtained from \cite{c3podata}, we extracted the point cloud. For each reconstruction in this subsection, the computational domain is defined as $\Omega = (0,0.5)\times(0,0.5)\times(0,1)$, discretized with a grid spacing of $h = 0.5/300$. Fig. \ref{c3po_volume} shows, respectively, the appearance of C-3PO in the movie \cite{c3pomovie}, the extracted point cloud, the reconstructed volume with $\epsilon = 2.4$e-$3$, and a cut-off view of the reconstructed volume.

\begin{figure}[htbp]
    \centering

    \begin{minipage}{1\linewidth}
        \centering
        \includegraphics[width=\linewidth]{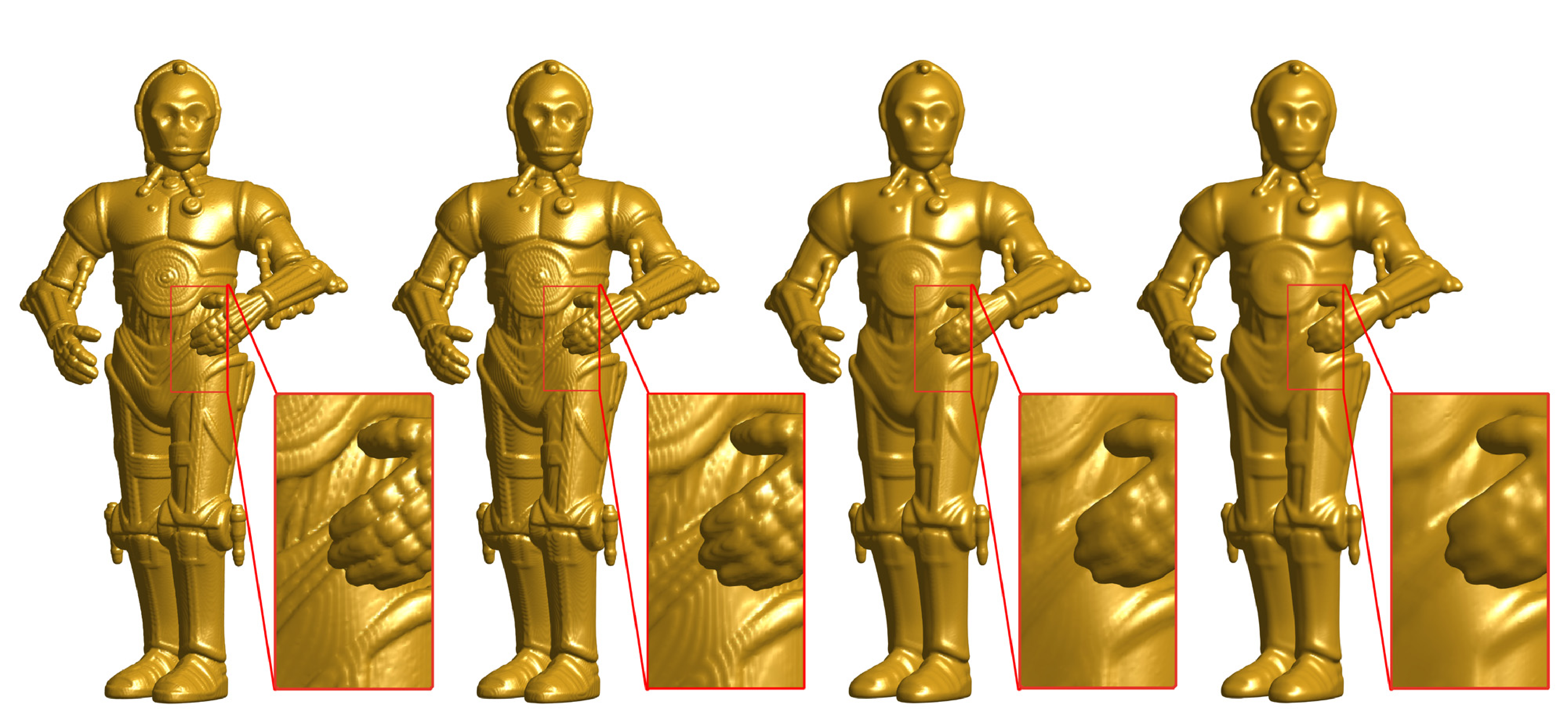}
    \end{minipage}

    \vspace{8mm}

    \begin{minipage}{1\linewidth}
        \centering
        \includegraphics[width=\linewidth]{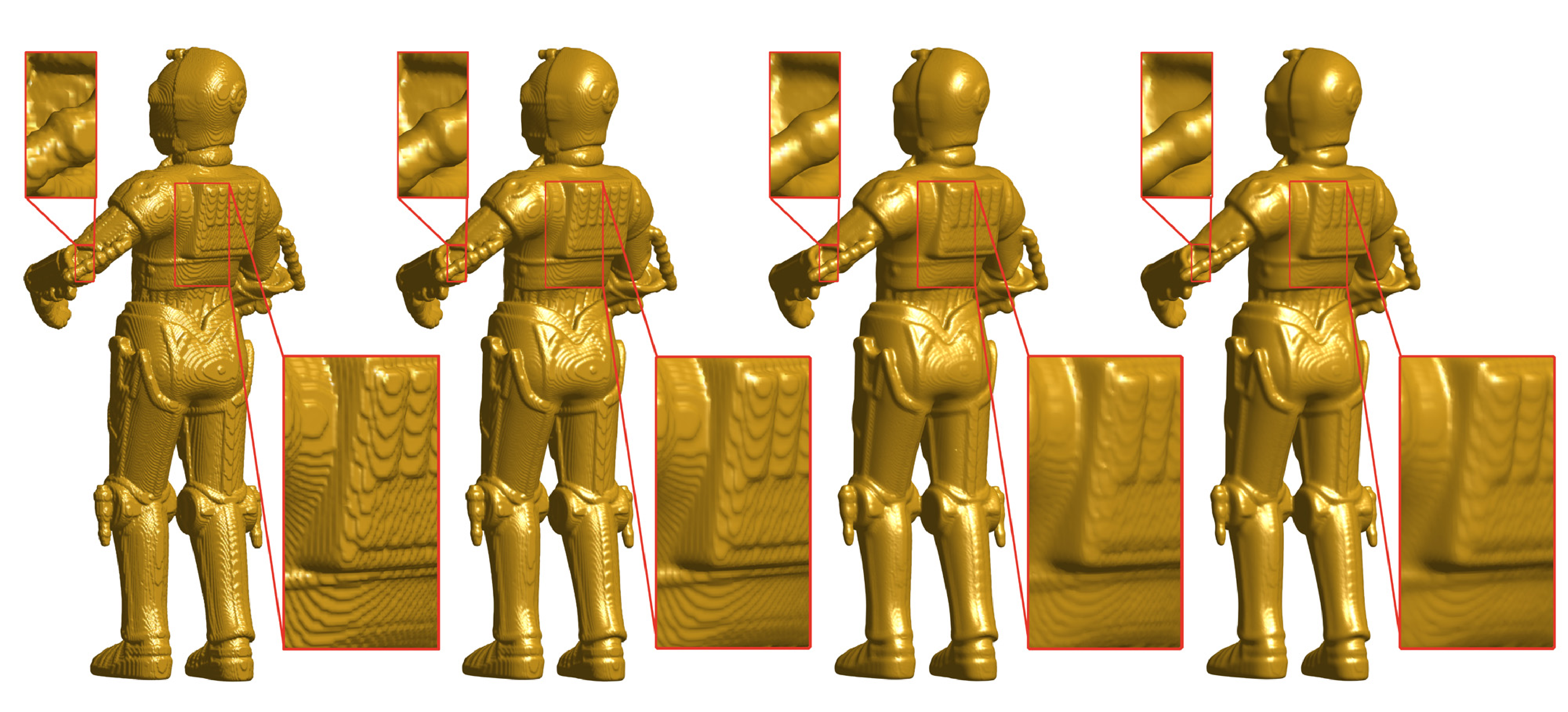}
        \par
        \vspace{0.1mm}
        \begin{minipage}{0.25\linewidth}
            \centering
        \end{minipage}
        \begin{minipage}{0.25\linewidth}
            \centering
        \end{minipage}
        \begin{minipage}{0.25\linewidth}
            \centering
        \end{minipage}
        \begin{minipage}{0.24\linewidth}
            \centering
        \end{minipage}
    \end{minipage}

    \caption{Impact of different $\epsilon$ values on the 3D reconstruction of C-3PO. From left to right, each column corresponds to the reconstructed volume of $\epsilon_a$, $\epsilon_b$, $\epsilon_c$, and $\epsilon_d$, respectively. The upper row displays the front views reconstructed with different $\epsilon$ values, including enlarged details of the hand region. The bottom row shows the back views under varying $\epsilon$ values, accompanied by enlarged views of the arm region and back.} \label{c3po_epsilon}
\end{figure}

To study the effects of different $\epsilon$ values on the reconstructed volume, each column in Fig. \ref{c3po_epsilon} presents the reconstruction results for $\epsilon_a$, $\epsilon_b$, $\epsilon_c$, and $\epsilon_d$, respectively. To facilitate visual comparison, the first row of Fig. \ref{c3po_epsilon} displays the front views of the reconstructed object, including enlarged details of the hands; the second row shows the back views, with enlarged views of the arm region and back.

By observing the experimental results, we find that as $\epsilon$ increases, the reconstructed model achieves better smoothness. In areas with intricate details but insufficient grid resolution, such as the hands in the first row of Fig. \ref{c3po_epsilon}, a larger $\epsilon$ may lead to loss of details due to increased smoothing. For the arm and back regions in the second row of Fig. \ref{c3po_epsilon}, which are initially not smooth, increasing $\epsilon$ from small values results in a smoother surface and more natural transitions. Therefore, for different models, we can empirically set the value of $\epsilon$. When reconstructing a model with many details, choosing a moderate $\epsilon$ is preferable to capture finer details. Conversely, if the initial reconstruction is relatively rough, selecting a larger $\epsilon$ can effectively smooth the surface and enhance reconstruction quality.

\subsection{Reconstruction of an actual object}
\begin{figure}[htbp]
    \begin{minipage}{0.5\linewidth}
    \centering
    \includegraphics[clip=true,width=2.81in]{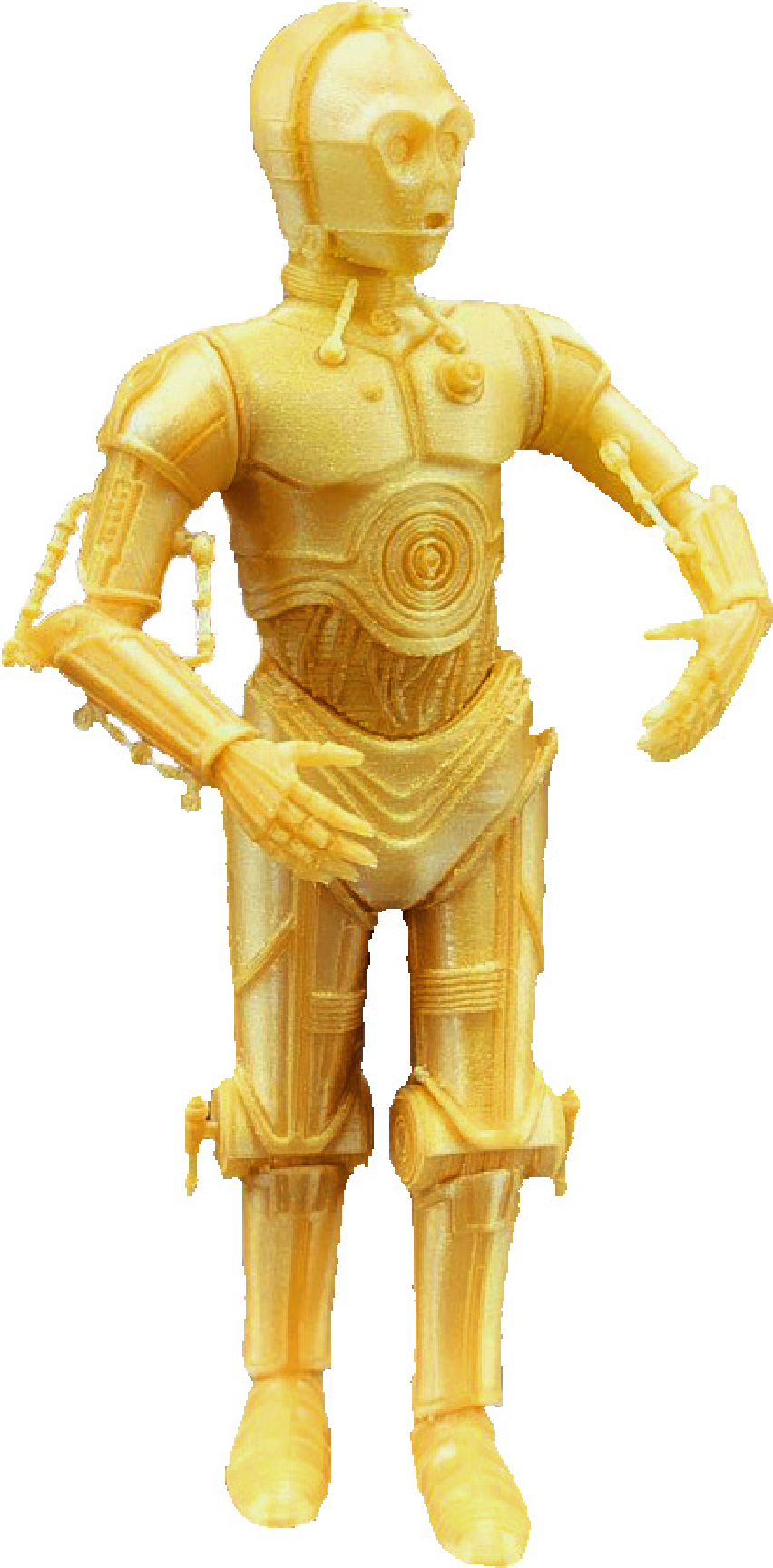}\\
    \vspace{3mm}
    \end{minipage}
    \begin{minipage}{0.5\linewidth}
    \centering
    \includegraphics[clip=true,width=2.71in]{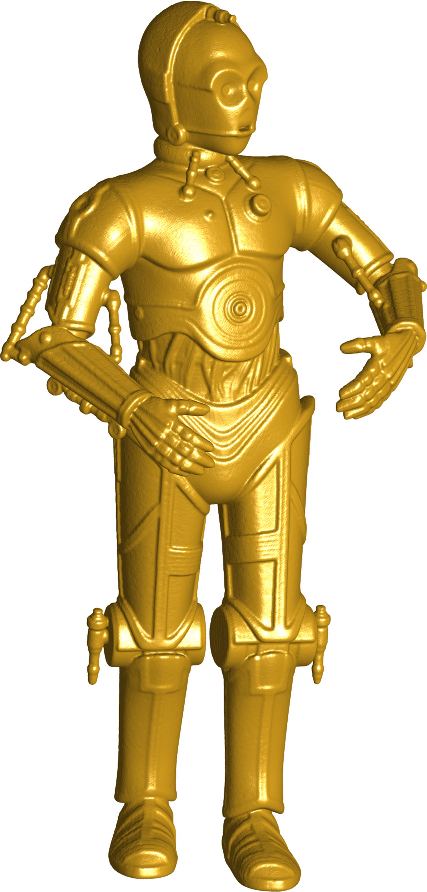}\\
    \end{minipage}\\
    \caption{Comparison with an actual 3D Object. The left column presents the actual object fabricated using a 3D printer, while the right column displays the numerical results of our proposed algorithm.} \label{c3poComparison}
    \end{figure}

In this subsection, we concisely compare our method with an actual C-3PO model created through 3D printing using the 3D file from \cite{c3podata}. C-3PO, as a humanoid robot, exhibits a surface with intricate topological features. The left column of Fig. \ref{c3poComparison} features an object we fabricated using a 3D printer. The right column displays the numerical outcome. To enhance the details of the reconstructed volume, we define the domain as $(0,0.5)\times(0,0.5)\times(0,1)$, discretized with a grid spacing of $h = 0.5/600$, and a transitional interface thickness parameter  $\epsilon = 5h/(2\sqrt{2}\tanh^{-1}(0.9))$. The actual and numerical results exhibit qualitative similarity.

\subsection{Impact of dispersed points on reconstruction}
\begin{figure}[htbp]
    \begin{minipage}{0.33\linewidth}
    \centering
    \includegraphics[clip=true,width=1.9in]{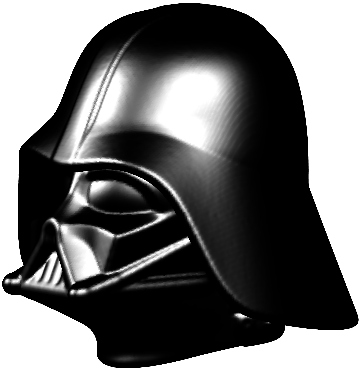}\\
    \vspace{2mm}
    \end{minipage}
    \begin{minipage}{0.33\linewidth}
    \centering
    \includegraphics[clip=true,width=1.9in]{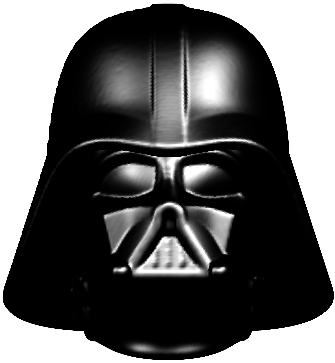}\\
    \vspace{2mm}
    \end{minipage}
    \begin{minipage}{0.32\linewidth}
    \centering
    \includegraphics[clip=true,width=2in]{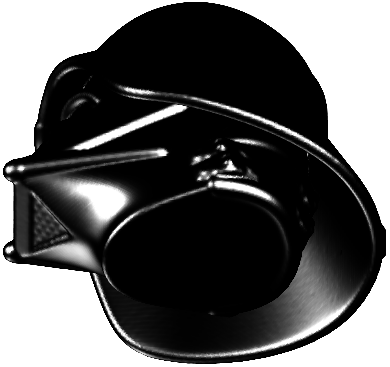 }\\
    \vspace{1.5mm}
    \end{minipage}\\
    \caption{Schematic illustration of Darth Vader's helmet. The reconstructed volumes are displayed from various perspectives, highlighting the algorithm's ability to capture intricate surface details.} \label{vader.eps}
    \end{figure}

\begin{figure}[htbp]
\begin{minipage}{0.33\linewidth}
\centering
\includegraphics[clip=true,width=1.9in]{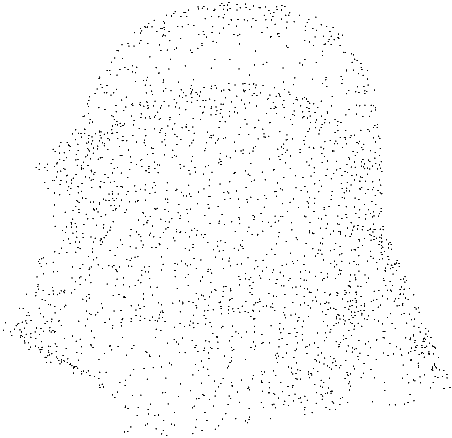}\\
\vspace{3mm}
\end{minipage}
\begin{minipage}{0.33\linewidth}
\centering
\includegraphics[clip=true,width=1.9in]{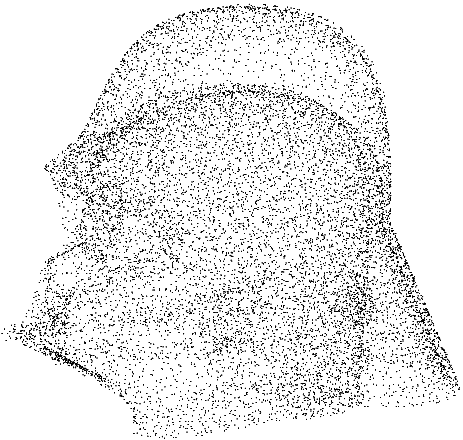}\\
\vspace{3mm}
\end{minipage}
\begin{minipage}{0.32\linewidth}
\centering
\includegraphics[clip=true,width=1.9in]{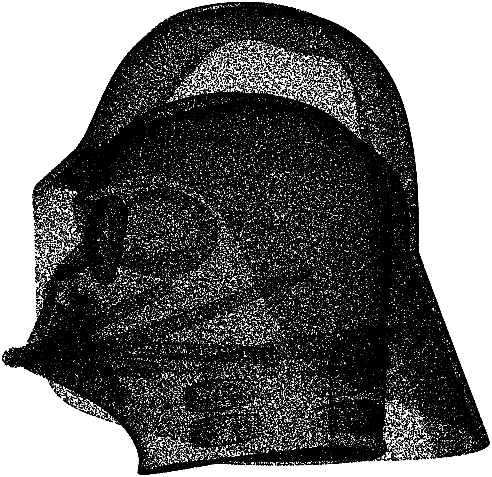}\\
\vspace{3mm}
\end{minipage}\\
\begin{minipage}{0.33\linewidth}
\centering
\includegraphics[clip=true,width=1.8in]{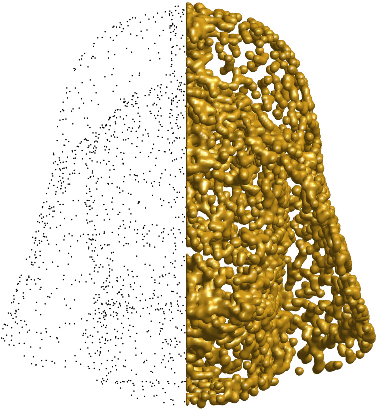}\\
\vspace{3mm}
\end{minipage}
\begin{minipage}{0.33\linewidth}
\centering
\includegraphics[clip=true,width=1.8in]{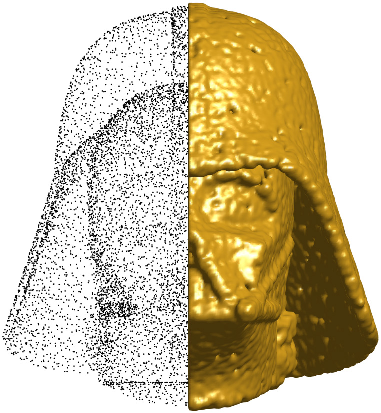}\\
\vspace{3mm}
\end{minipage}
\begin{minipage}{0.32\linewidth}
\centering
\includegraphics[clip=true,width=1.8in]{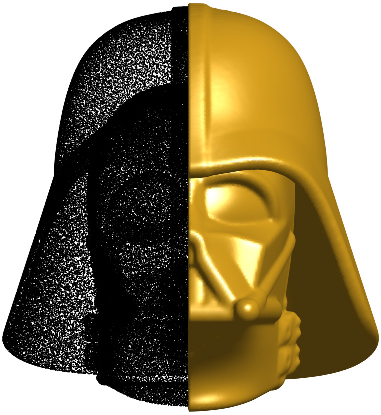}\\
\vspace{3mm}
\end{minipage}\\
\begin{minipage}{0.33\linewidth}
\centering
\includegraphics[clip=true,width=1.8in]{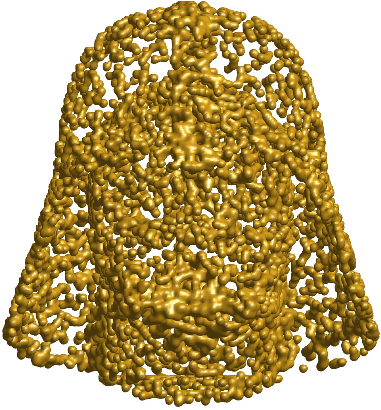}\\
$5000$
\end{minipage}
\begin{minipage}{0.33\linewidth}
\centering
\includegraphics[clip=true,width=1.8in]{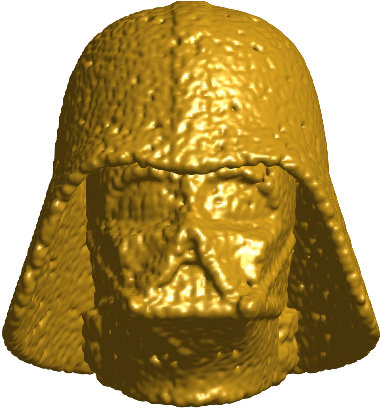}\\
$25000$
\end{minipage}
\begin{minipage}{0.32\linewidth}
\centering
\includegraphics[clip=true,width=1.8in]{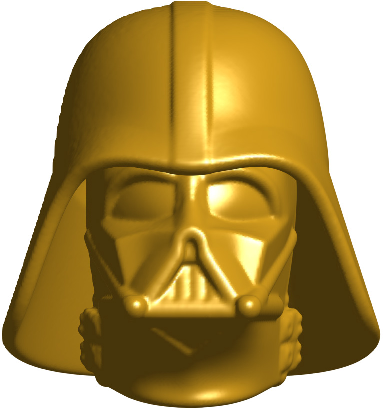}\\
$500000$
\end{minipage}
\caption{Impact of data point quantity on three-dimensional volume reconstruction. The upper row display the side-view datasets with differing numbers of points. The middle row present cut-off views of the associated 3D profiles from the front perspective. The bottom row illustrates the complete 3D profiles from the front view. The numeber of data points is specified beneath each figure. } \label{turtt1.eps}
\end{figure}

We use the Darth Vader's helmet, sourced from \cite{Vaderdata1} to show the effect of points cloud on reconstruction. After preprocessing the original 3D file, we obtained a point cloud suitable for reconstruction, consisting of $500000$ data points. The domain is defined as $(0,1)\times(0,1)\times(0,1)$, with parameters $\hat{\delta t} = 2.5$e-$5$, $h = 0.1/40$, and $\epsilon = 0.0036$. Darth Vader is an iconic character from the renowned movie series \textit{Star Wars}. To show the capability of our proposed scheme, we give a schematic illustration in Fig. \ref{vader.eps}, where the helmet is reconstructed using our method. To more accurately depict Darth Vader's appearance, we colored the reconstructed volume black using MATLAB. By reducing the number of points by factors of $20$ and $100$, we obtained two additional sets of coarser point clouds. The top row of Fig. \ref{turtt1.eps} shows the side-view point cloud data sets with decreasing numbers of points. For enhanced visualization, only a portion of the points is displayed. The middle row of Fig. \ref{turtt1.eps} presents the corresponding front cut-off views of the reconstructed 3D profiles for the different data point counts. With the reduction of the number of points, it is evident that the resulting 3D profiles are less smooth. Consequently, to achieve precise results, the present algorithm requires a sufficient number of scattered points.

\subsection{Sensitivity of $\epsilon$ to point density}
\begin{figure}[htbp]
    \centering

    \begin{minipage}{0.04\linewidth}\raggedleft
    {\Large$\epsilon_1$}
    \end{minipage}%
    \begin{minipage}{0.3\linewidth}\centering
    \includegraphics[clip=true,width=1.2in]{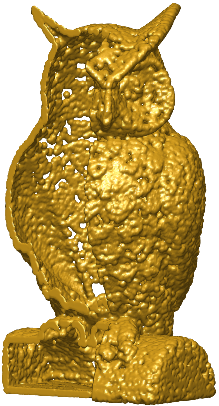}\\
    \vspace{1mm}
    \end{minipage}%
    \begin{minipage}{0.3\linewidth}\centering
    \includegraphics[clip=true,width=1.2in]{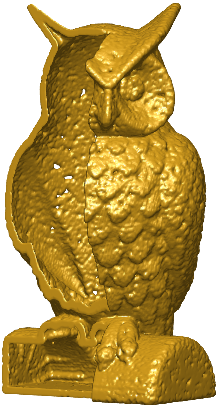}\\
    \vspace{1mm}
    \end{minipage}%
    \begin{minipage}{0.3\linewidth}\centering
    \includegraphics[clip=true,width=1.2in]{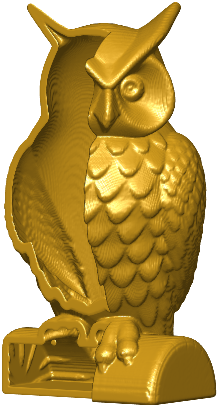}\\
    \vspace{1mm}
    \end{minipage}\\[1mm]

    \begin{minipage}{0.04\linewidth}\raggedleft
    {\Large$\epsilon_2$}
    \end{minipage}%
    \begin{minipage}{0.3\linewidth}\centering
    \includegraphics[clip=true,width=1.2in]{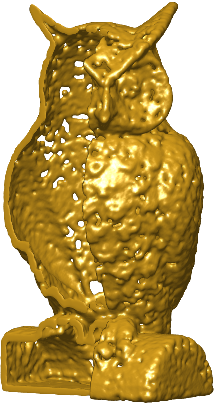}\\
    \vspace{1mm}
    \end{minipage}%
    \begin{minipage}{0.3\linewidth}\centering
    \includegraphics[clip=true,width=1.2in]{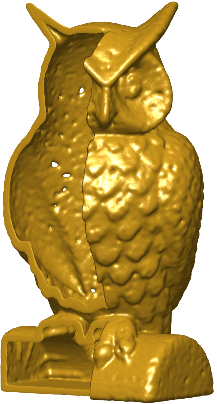}\\
    \vspace{1mm}
    \end{minipage}%
    \begin{minipage}{0.3\linewidth}\centering
    \includegraphics[clip=true,width=1.2in]{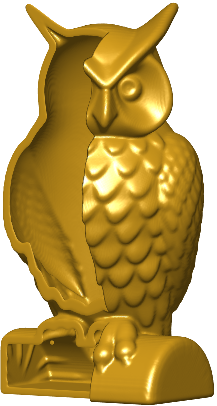}\\
    \vspace{1mm}
    \end{minipage}\\[1mm]

    \begin{minipage}{0.04\linewidth}\raggedleft
    {\Large$\epsilon_3$}
    \end{minipage}%
    \begin{minipage}{0.3\linewidth}\centering
    \includegraphics[clip=true,width=1.2in]{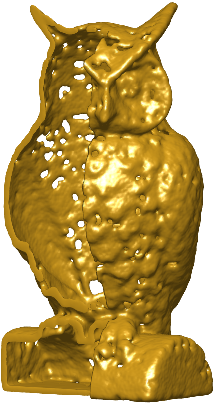}\\
    $17100$
    \end{minipage}%
    \begin{minipage}{0.3\linewidth}\centering
    \includegraphics[clip=true,width=1.2in]{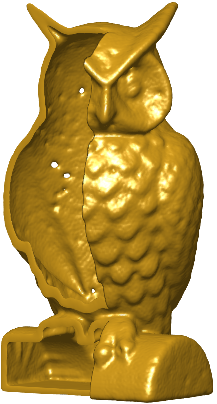}\\
    $30780$
    \end{minipage}%
    \begin{minipage}{0.3\linewidth}\centering
    \includegraphics[clip=true,width=1.2in]{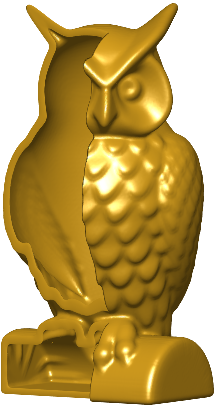}\\
    $153900$
    \end{minipage}\\

    \caption{Impact of point density and interface width $\epsilon$ on 3D reconstruction. The upper, middle, and bottom rows correspond to $\epsilon_1 = 9$e-$4$, $\epsilon_2 = 1.8$e-$3$, and $\epsilon_3 = 3.6$e-$3$, respectively; the left, center, and right columns correspond to different point densities. For clearer visualization, one quarter of each reconstruction is clipped. The number of data points is shown beneath each column.}
    \label{fig:nine_comparison_new}
    \end{figure}

In this subsection, we conduct a $3\times 3$ two-factor study to provide practical guidance for choosing the interface thickness $\epsilon$ under different point densities. We use the \emph{owl} point cloud as the testbed.
The computational domain is $\Omega=(0,0.6)\times(0,0.6)\times(0,1.0)$, discretized on a uniform grid $(N_x,N_y,N_z)=(240,240,400)$, which gives the mesh size $h=0.6/240=2.5\times 10^{-3}$.
Figure~\ref{fig:nine_comparison_new} arranges the reconstructions so that the columns correspond to three points densities ($17100$, $30780$, $153900$), while rows correspond to three interface widths $\epsilon_1 = 1.5h/(2\sqrt{2}\tanh^{-1}(0.9)) = 9$e-$4$, $\epsilon_2 = 3h/(2\sqrt{2}\tanh^{-1}(0.9)) = 1.8$e-$3$, and $\epsilon_3 = 6h/(2\sqrt{2}\tanh^{-1}(0.9)) = 3.6$e-$3$.
For clearer visualization, one quarter of each reconstructed volume is clipped away. From the results we observe that, when the point sampling is sparse, increasing $\epsilon$ yields visibly cleaner and more coherent surfaces. This behavior is consistent with the role of a larger diffuse-interface thickness in attenuating the roughness and small voids introduced by undersampling in the point cloud.

    \subsection{Volume evolution under varying $\epsilon$}
\begin{figure}[htbp]
        \centering
        \includegraphics[clip=true,width=0.95\linewidth]{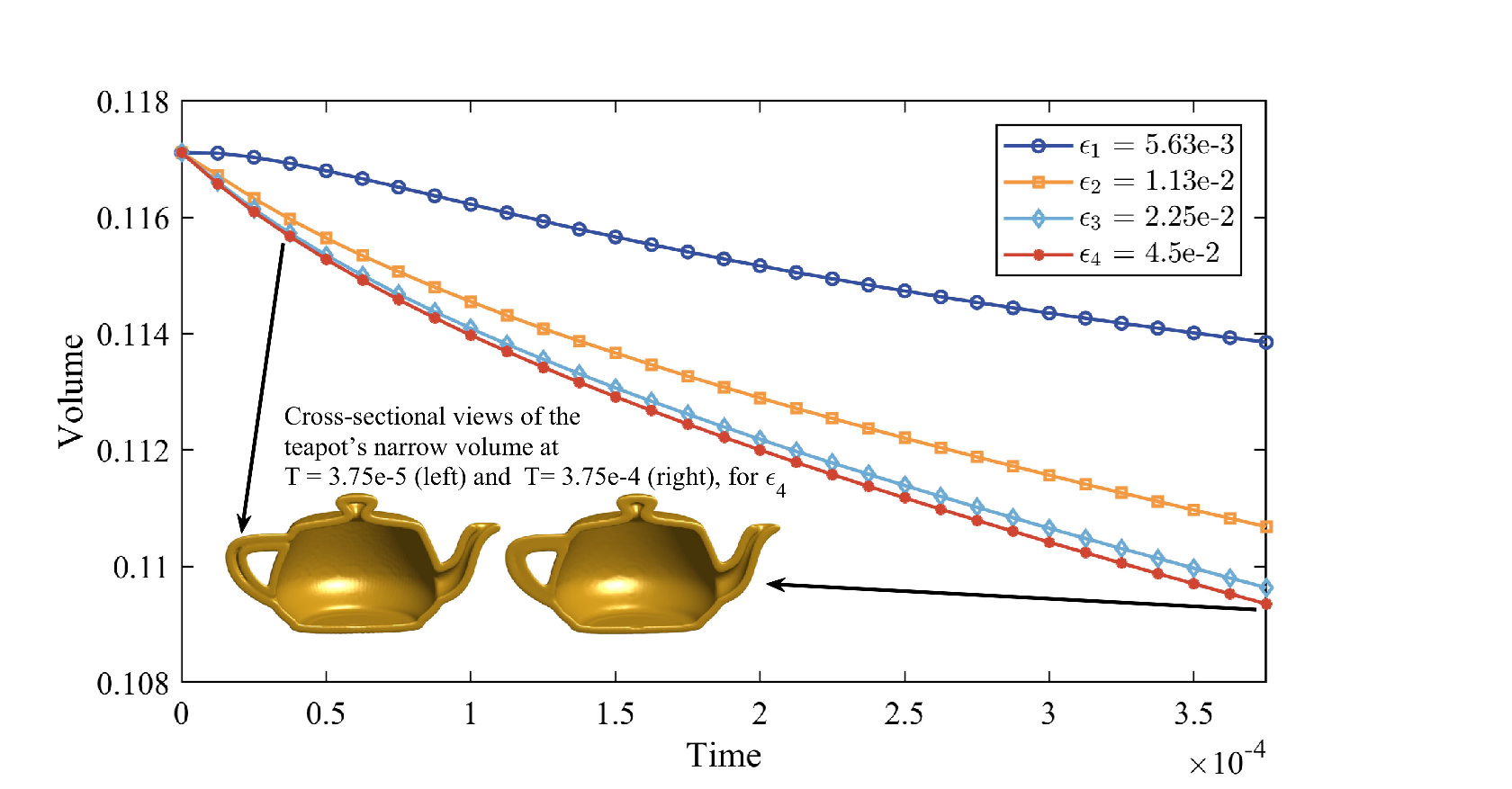}
        \caption{Volume evolution under varying $\epsilon$.}
        \label{vol-eps-comparison}
    \end{figure}

    In this subsection, we examine the evolution of the diffuse volume of the \emph{teapot} geometry for four interface-thickness parameters $\epsilon$ that increase by factors of two. The volume is defined via the phase-field variable $\phi(\mathbf{x},t)\in[-1,1]$ as the integral of the diffuse indicator $(1+\phi)/2$ over $\Omega$,
    \begin{equation}\label{eq:vol_def_new}
      V(t) \;=\; \int_{\Omega} \frac{1+\phi(\mathbf{x},t)}{2}\,\mathrm{d}\mathbf{x},
    \end{equation}
    and, on a uniform Cartesian mesh of spacing $h$, the discrete value at time level $t_n$ is $V^n=h^3\sum_{i,j,k}\bigl(1+\phi^n_{i,j,k}\bigr)/2$. We advance with a constant time step $\hat{\delta t}=1.25$e-$5$ and record the volume at 30 successive steps, up to $T=3.75$e-$4$. We use four successively doubled interface thicknesses, $\epsilon_1=5.63$e-$3$, $\epsilon_2=1.13$e-$2$, $\epsilon_3=2.25$e-$2$, and $\epsilon_4=4.5$e-$2$, while all other teapot parameters follow the earlier definition. Fig.~\ref{vol-eps-comparison} reports the corresponding curves $V(t)$: the volume exhibits a monotone decrease for all $\epsilon$; for each $t$ we observe $V_{\epsilon_1}(t)\ge V_{\epsilon_2}(t)\ge V_{\epsilon_3}(t)\ge V_{\epsilon_4}(t)$. Larger $\epsilon$ produces stronger curvature-driven smoothing and thus a faster loss of volume in narrow features, whereas the gaps between the curves diminish as $\epsilon$ decreases, indicating reduced $\epsilon$-sensitivity and consistency with convergence to the sharp-interface limit. For illustration, two cross-sectional views of the teapot's narrow region for $\epsilon_4$ at $T=3.75$e-$5$ (left) and $T=3.75$e-$4$ (right) are included in the figure, showing progressive neck-thinning aligned with the measured decline in $V(t)$.

    \subsection{Effect of the stabilization constant $S$}\label{subsec:S_effect}
\begin{figure}[htbp]
    \centering

    \begin{minipage}{0.32\linewidth}
        \centering
        \includegraphics[clip=true,width=1.8in]{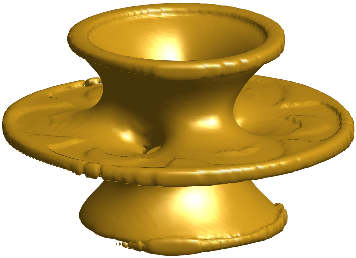}
        \\$S_1 = \frac{0}{\epsilon^2}$\vspace{2mm}
    \end{minipage}
    \begin{minipage}{0.32\linewidth}
        \centering
        \includegraphics[clip=true,width=1.8in]{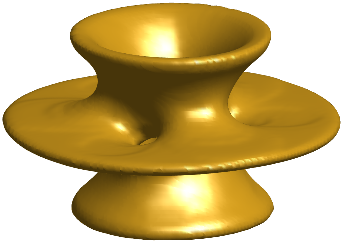}
        \\$S_2 = \frac{2}{\epsilon^2}$\vspace{2mm}
    \end{minipage}
    \begin{minipage}{0.32\linewidth}
        \centering
        \includegraphics[clip=true,width=1.8in]{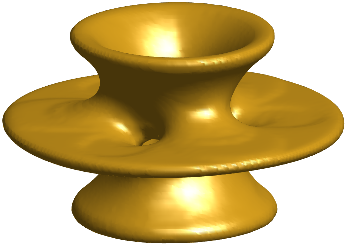}
        \\$S_3 = \frac{4}{\epsilon^2}$\vspace{2mm}
    \end{minipage}
    \\[-1mm]

    \caption{Effect of the stabilization constant $S$ on the Costa reconstruction with a large time step $\hat{\delta t} = 1.25$e-$3$. From left to right: $S_1=\frac{0}{\epsilon^2}$ (no stabilization; unstable and strongly degraded), $S_2=\frac{2}{\epsilon^2}$, and $S_3=\frac{4}{\epsilon^2}$. Both stabilized choices are visually stable; the larger $S$ yields slightly more diffusion.}
    \label{diff_S}
\end{figure}

In this subsection, we assess the influence of the stabilization constant $S$ on stability and accuracy of the reconstruction. We use the Costa--Hoffman--Meeks surface as the test case and keep all parameters the same as in the preceding subsection, except that we deliberately set a relatively large time step $\hat{\delta t} = 1.25$e-$3$ to stress the stability of the update. Figure~\ref{diff_S} reports the reconstructions for three choices: $S_1 = \frac{0}{\epsilon^2}$ (i.e., no stabilization), $S_2 = \frac{2}{\epsilon^2}$, and $S_3 = \frac{4}{\epsilon^2}$. When $S=0$, the computation becomes unstable under this large $\hat{\delta t}$, producing severe artifacts. In particular, the discrete energy fails to remain bounded and the iteration diverges under this step size. In contrast, both $S=\frac{2}{\epsilon^2}$ and $S=\frac{4}{\epsilon^2}$ yield visually stable results.

    The $S$-term acts as a linear damping that enlarges the scheme's stability region and permits larger time steps. However, it also introduces a consistency error in time tied to the second-order temporal residual (proportional to $\partial_{tt}\phi$); an excessively large $S$ may over-smooth fine features and degrade accuracy. In practice, we found $S=\frac{2}{\epsilon^2}$ to be a robust compromise: it stabilizes the CN update for the time steps while avoiding noticeable loss of detail. Using $S=\frac{4}{\epsilon^2}$ provides extra robustness for more aggressive $\hat{\delta t}$, at the cost of more diffusion. Conversely, setting $S=0$ is viable only when $\hat{\delta t}$ is chosen sufficiently small. These results suggest that $S = \frac{2}{\epsilon^2}$ is empirically a good choice.

\subsection{Comparison with SAV scheme}
\begin{figure}[htbp]
    \begin{minipage}{0.33\linewidth}
        \centering
        \includegraphics[clip=true,width=1.9in]{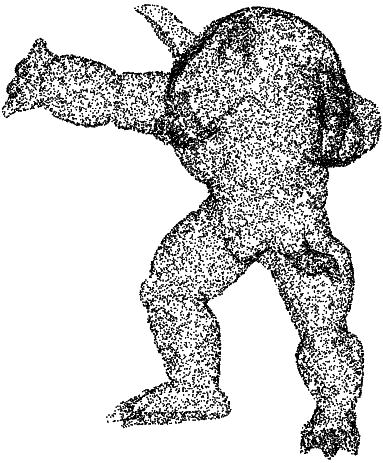}\\
        \vspace{2mm}
        \end{minipage}
    \begin{minipage}{0.33\linewidth}
    \centering
    \includegraphics[clip=true,width=1.9in]{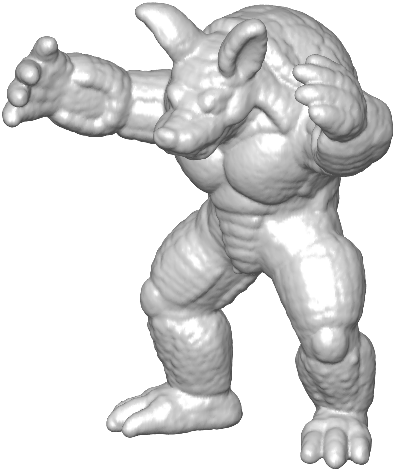}\\
    \vspace{2mm}
    \end{minipage}
    \begin{minipage}{0.33\linewidth}
    \centering
    \includegraphics[clip=true,width=1.9in]{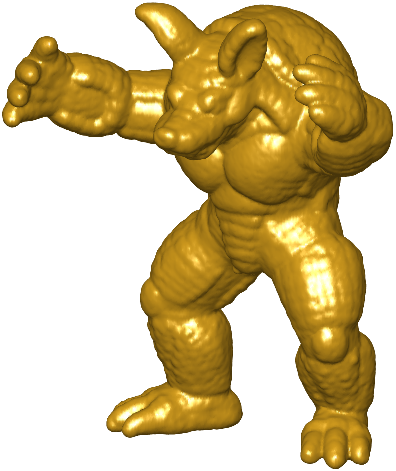}\\
    \vspace{2mm}
    \end{minipage}\\
    \caption{3D reconstruction of the Armadillo model. From left to right: points cloud, reconstruction using SAV scheme, reconstruction using Lagrange-multiplier scheme.} \label{sav-lm-reconstruction}
    \end{figure}

\begin{figure}[htbp]
    \centering
    \includegraphics[clip=true,width=0.95\linewidth]{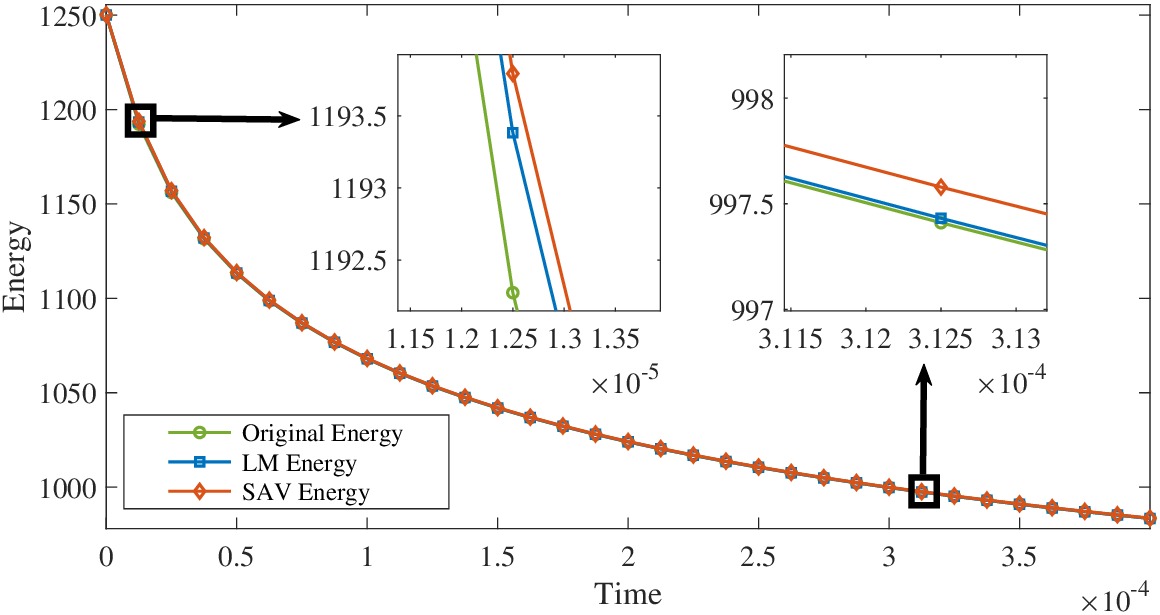}
    \caption{Energy dissipation of the original, Lagrange multiplier, and SAV schemes.}
    \label{sav-comparison}
\end{figure}

In this subsection, we implement a comparison with the SAV scheme \cite{SAVV1,SAVV2}. As a variant of SAV scheme, the Lagrange multiplier (LM) scheme not only inherits the merit of energy stability resulting from SAV scheme, but also increase the consistency between original and numerical energy. To show this advantage, we use the points cloud of Armadillo and generate the reconstructions with SAV scheme and LM scheme, respectively. Figure \ref{sav-lm-reconstruction} shows the points cloud and reconstruction results, the 3D profiles are highly similar because we use a relatively fine time step $\Delta t = 1.25$e-$5$. In Fig. \ref{sav-comparison}, the evolutions of original energy and numerical energy obtained by SAV and LM schemes are plotted. From the close-up views, we can observe that the energy consistency is obviously improved by the LM scheme.

\subsection{Comparison with prior method for 3D reconstruction}\label{comparison_with_BDF2}
In this subsection, we compare the proposed method with the energy-stable BDF2 strategy~\cite{Cai2025}. Both approaches are formulated on the distance–weighted Allen–Cahn model for narrow-volume reconstruction and use finite differences in space.

\begin{figure}[htbp]
    \centering

    \begin{minipage}{0.30\linewidth}
        \vspace{5mm}
        \centering
        \includegraphics[clip=true,width=2.0in]{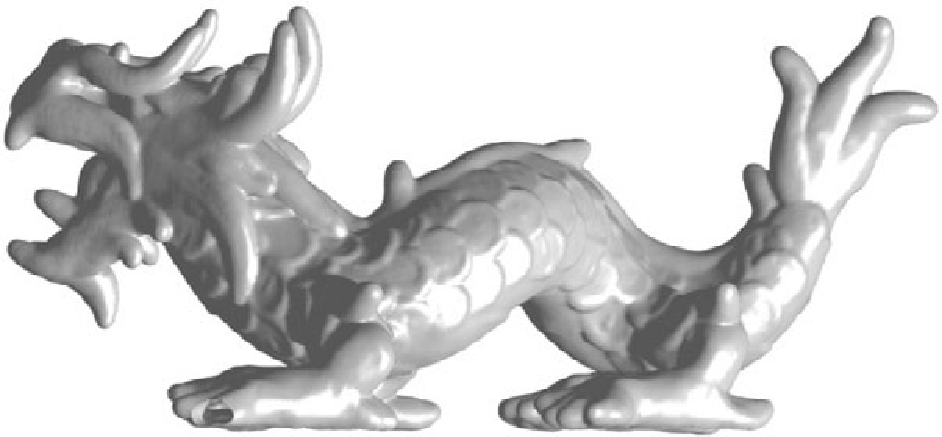}
        \vspace{2mm}
    \end{minipage}
    \begin{minipage}{0.30\linewidth}
        \centering
        \includegraphics[clip=true,width=1.6in]{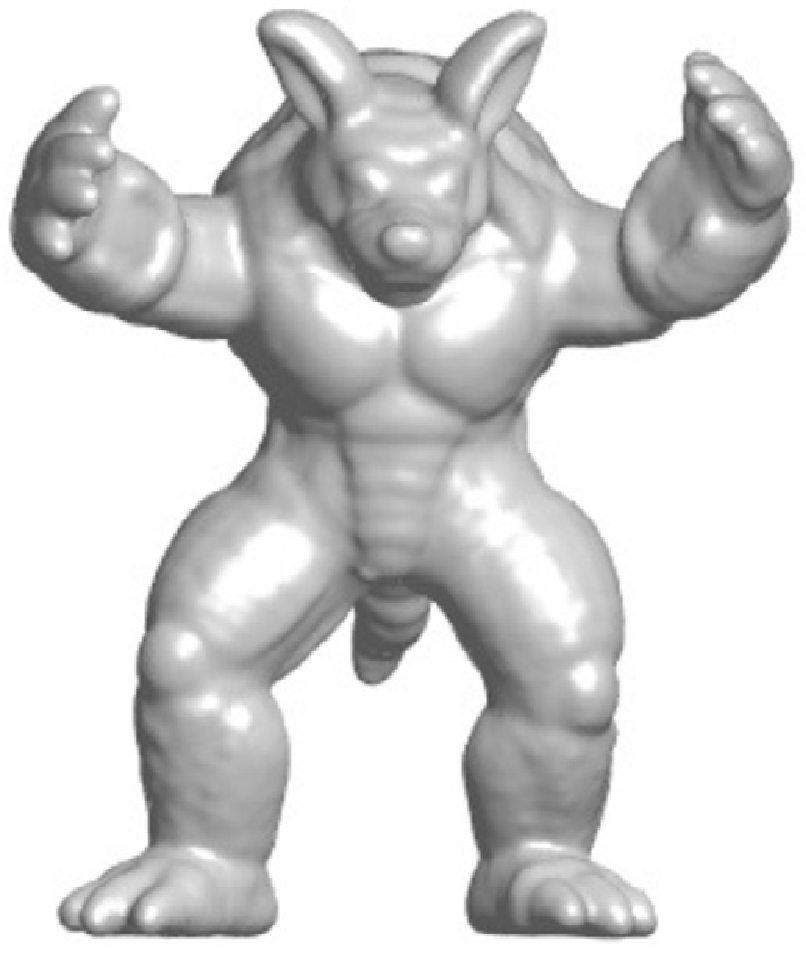}
        \vspace{2mm}
    \end{minipage}
    \begin{minipage}{0.30\linewidth}
        \vspace{5mm}
        \centering
        \includegraphics[clip=true,width=2.0in]{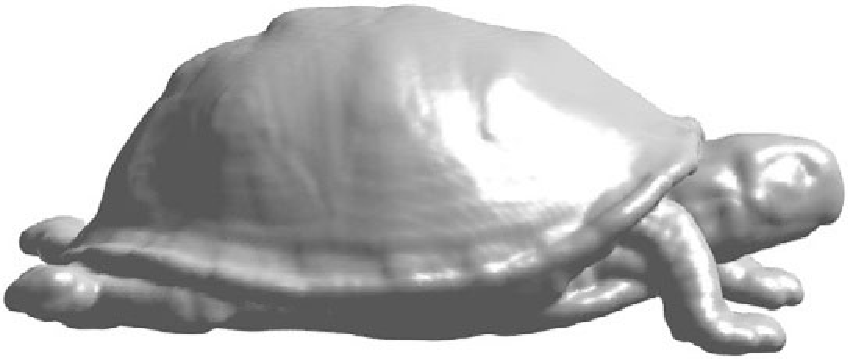}
        \vspace{2mm}
    \end{minipage}
    \\[1mm]

    \begin{minipage}{0.30\linewidth}
        \vspace{5mm}
        \centering
        \includegraphics[clip=true,width=2.0in]{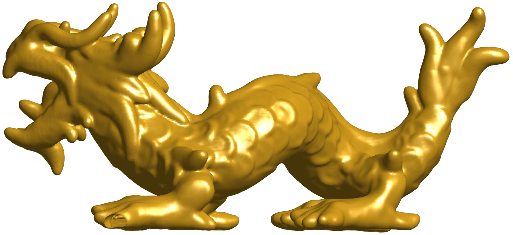}
        \vspace{2mm}
    \end{minipage}
    \begin{minipage}{0.30\linewidth}
        \centering
        \includegraphics[clip=true,width=1.6in]{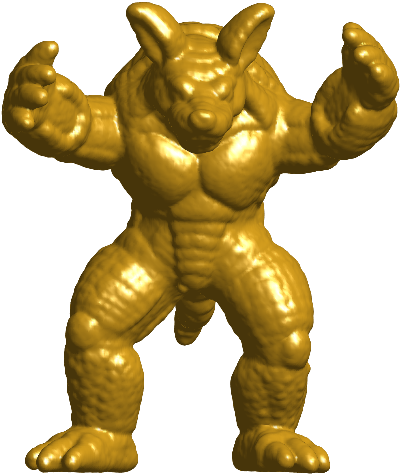}
        \vspace{2mm}
    \end{minipage}
    \begin{minipage}{0.30\linewidth}
        \vspace{5mm}
        \centering
        \includegraphics[clip=true,width=2.0in]{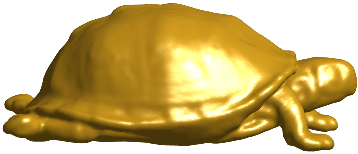}
        \vspace{2mm}
    \end{minipage}

    \caption{Three-dimensional narrow-volume reconstructions of the Asian Dragon, Armadillo, and Turtle models for comparison. The upper row shows results taken from ~\cite{Cai2025} with permission from the corresponding author, whereas the lower row displays the outcomes generated in this study.}
    \label{comparisonBDF2}
\end{figure}

For a fair assessment, we adopt the same point clouds and numerical parameters as in~\cite{Cai2025} for the Asian Dragon, Armadillo, and Turtle models. Figure~\ref{comparisonBDF2} places the results side-by-side: the \emph{top} row shows the BDF2 reconstructions reported by Cai \emph{et al.}~\cite{Cai2025}, while the \emph{bottom} row displays our reconstructions. Under identical settings, the two sets of surfaces are highly consistent, and visually our results preserve slightly more fine-scale details.

Algorithmically, our Lagrange--multiplier CN discretization differs in several structure-preserving aspects that are advantageous in practice: (i) it enforces the discrete energy law of the \emph{original} AC functional via a time-dependent Lagrange multiplier, avoiding the extrapolation/truncation used to linearize the cubic term in the BDF2 scheme; (ii) each step decouples into two linear subproblems plus a single scalar constraint, making the method self-starting and free of auxiliary stabilization parameters. These features yield unconditional stability with simple linear solves and minimize parameter-tuning overhead in implementation.

\subsection{Algorithmic complexity and convergence criteria}
To quantify algorithmic complexity and convergence behavior, we select the \textbf{C-3PO} dataset as a representative test case. The computational domain is $\Omega=(0,0.5)\times(0,0.5)\times(0,1)$, discretized on a uniform grid $(N_x,N_y,N_z)=(300,300,600)$ with spacing $h=0.5/300$ (so $h=l_x/N_x=l_y/N_y=l_z/N_z$). We take the time step $\hat{\delta t} = 6.25$e-$6$ and integrate for $T=10\,\hat{\delta t}$ (10 steps). Here, the two linear substeps correspond to the affine-in-$Q$ decoupling: the ``$u$'' subproblem solves Eqs.~\eqref{eq:update_phi1}--\eqref{eq:update_mu1} (or equivalently the operator form \eqref{0CN1}) for $(\phi_1,\mu_1)$, and the ``$v$'' subproblem solves Eqs.~\eqref{eq:update_phi2}--\eqref{eq:update_mu2} (or \eqref{0CN2}) for $(\phi_2,\mu_2)$. Under this setting, the GS iterations for the two substeps average $16.00$ (``$u$'') and $15.00$ (``$v$''), i.e., $31.00$ total GS iterations per time step on average, while the Newton solver for the scalar Lagrange multiplier $Q^{n+\frac{1}{2}}$ requires on average $18.10$ iterations per update.

\begin{table}[htbp]
\centering
\caption{GS iterations and Newton iterations over 10 time steps for the C-3PO test.}
\begin{tabular}{lccccccccccc}
\toprule
 & 1 & 2 & 3 & 4 & 5 & 6 & 7 & 8 & 9 & 10 & Avg \\
\midrule
GS iterations (u) & 20 & 20 & 15 & 15 & 15 & 15 & 15 & 15 & 15 & 15 & 16.00 \\
GS iterations (v) & 15 & 15 & 15 & 15 & 15 & 15 & 15 & 15 & 15 & 15 & 15.00 \\
Newton iterations & 1 & 1 & 172 & 1 & 1 & 1 & 1 & 1 & 1 & 1 & 18.10 \\
Total GS iterations & 35 & 35 & 30 & 30 & 30 & 30 & 30 & 30 & 30 & 30 & 31.00 \\
\bottomrule
\end{tabular}
\label{tab:gs_newton_c3po}
\end{table}

As summarized by the \emph{Avg} column in Table~\ref{tab:gs_newton_c3po}, the averages are $16.00$ (GS iterations in ``$u$''), $15.00$ (GS iterations in ``$v$''), $31.00$ (total GS iterations per step), and $18.10$ (Newton iterations per update). We note that the unusually large Newton count at step 3 ($172$ iterations) arises from the early-time transient when enforcing the discrete energy constraint \eqref{cnQdf} via the polynomial relation \eqref{simplified}: the Jacobian with respect to $Q^{n+\frac{1}{2}}$ becomes poorly conditioned (the energy residual is highly nonlinear and its derivative is small near the initial profile), so an undamped Newton with the default initial guess requires many corrections. After the transient, the solution is close to the descent manifold and $Q^{n+\frac{1}{2}}$ stays near $1$, leading to rapid (often single-iteration) updates in subsequent steps. These results collectively indicate robust convergence and practical efficiency of the decoupled CN scheme on this high-resolution setup.

\subsection{CPU time for representative reconstructions}

\begin{table}[htbp]
\centering
\caption{CPU time statistics for six points cloud reconstructions (seconds).}
\label{tab:cpu_time}
\begin{tabular}{lcccc}
\toprule
Points cloud & Points & Grid size & CPU time (IPF) [s] & Avg CPU time (PTS) [s] \\
\midrule
Teapot &26,103 & 4,194,304  & 0.43 & 2.50 \\
Costa & 184,478 & 2,097,152  & 2.69 & 1.27 \\
Horse & 40,116 & 4,194,304  & 0.65 & 2.61 \\
happy Buddha & 540,616 & 33,554,432  & 8.73 & 27.66 \\
Armadillo & 129,732 & 11,000,000  & 2.04 & 6.33 \\
Stanford Dragon & 218,823 & 30,720,000  & 3.56 & 26.15 \\
\bottomrule
\end{tabular}
\end{table}
To quantify computational cost, we report wall-clock CPU times for six representative point-cloud reconstructions (Teapot, Costa--Hoffman--Meeks, Horse, happy Buddha, Armadillo, Stanford Dragon). All experiments were executed on a laptop (Lenovo Yoga Pro 14s, 2023) equipped with an Intel Core i9-13905H CPU, without GPU acceleration. For each dataset we provide (i) the time to construct the initial phase field (IPF), which includes the unsigned-distance-based initialization and edge weight $g(\mathbf{x})$, and (ii) the average per-time-step (PTS) cost, obtained by averaging the wall-clock time of the first five steps (i.e., $5\,\Delta t$). Abbreviations and units are indicated in Table~\ref{tab:cpu_time} (IPF [s], PTS [s]).

The IPF time increases with both the grid size and the number of points, because the distance-based initialization is evaluated over the computational grid and depends on the input point cloud. By contrast, the PTS cost is largely insensitive to the number of points (the point cloud is not queried during time stepping) and is primarily governed by the grid resolution; consequently, the PTS cost grows approximately linearly with the number of cells in the computational domain. In Table~\ref{tab:cpu_time}, the Grid size column reports $N_x\times N_y\times N_z$ (number of cells in the computational domain). Consistent with these observations, the IPF time increases as either the grid is refined or the point sampling becomes denser, while the PTS time increases primarily with grid refinement. These measurements reflect the practical costs on commodity hardware and demonstrate the feasibility and efficiency of the proposed solver for moderate-to-large point-cloud reconstructions.

\section{Conclusions and discussions}\label{sec5}
We developed a practical and stable method incorporating a Lagrange multiplier approach for 3D reconstruction. The equilibrium equations were constructed by introducing a time-dependent variable $Q$. A stabilized CN-type method was considered to update the solutions in time. The spatial discretization was described in detail. The analytical proof of energy stability for the fully discrete scheme was implemented. Several benchmarks were used to validate the unconditional stability and the desired accuracy. Reconstructions of various complex objects demonstrated the practicality of our algorithm. Moreover, we conducted tests with various values of transitional interface thickness parameters using the point cloud of C-3PO, exploring the effects of these parameters on the reconstruction. Subsequently, we compared our numerical results with an actual 3D object, revealing qualitative similarities between them. Furthermore, we explored the effect of the number of points on reconstruction using the point cloud of Darth Vader's helmet, demonstrating that a sufficient number of data is essential for accurate reconstruction. In addition, we investigated the sensitivity of the reconstruction to the diffuse-interface thickness $\epsilon$ and to point density, and tracked the volume evolution under varying $\epsilon$, which together provide practical guidelines for parameter selection. We further examined the role of the stabilization constant $S$, observing a broad stability window and identifying $S=2/\epsilon^{2}$ as a robust default. Comparative studies against an SAV-based scheme and a prior BDF2-type method (see Subsection~\ref{comparison_with_BDF2}) confirmed the advantages of the proposed decoupled CN approach in terms of energy dissipation, accuracy, and practical efficiency. We also summarized the algorithmic complexity and convergence criteria, and reported CPU timings on commodity hardware, showing that the per-step cost scales mainly with grid resolution and is largely insensitive to the number of input points. Taken together, these results indicate that the proposed framework is practical, reliable, and ready for robust 3D reconstruction from unorganized point clouds.

It is worth noting that the performance of present reconstruction method depends on the quality of points cloud. If there are obvious noises in the points cloud, the reconstructing 3D profile might be rough or defective. In our upcoming work, a 3D smoothing technique will be considered to reduce the effect of noises.

\textit{Future work.} In future research, especially when targeting large-scale three-dimensional reconstructions or when preparing for engineering deployment, we will adopt adaptive mesh refinement (AMR), whose multi-resolution discretization allows local refinement near complex geometry while controlling memory and computational costs. We will integrate an AMR-based volumetric representation with scalable data structures and parallel I/O to support computation and storage at scale. This is expected to improve efficiency while preserving the method's accuracy and stability.

\section*{Acknowledgment}
Junxiang Yang is supported by Macau Science and Technology Development Fund (FDCT) (No. FDCT-24-080-SCSE). The authors appreciate the anonymous reviewers for their constructive comments.


\begin{thebibliography}{99}

\bibitem{shapeRe01}
S.S. Jinka, A. Srivastava, C. Pokhariya, A. Sharma, P.J. Narayanan,
SHARP: Shape-Aware Reconstruction of People in Loose Clothing,
International Journal of Computer Vision 131 (2023) 918--937.


\bibitem{shapeRe02}
V. Leroy, J.-S. Franco, E. Boyer,
Volume Sweeping: Learning Photoconsistency for Multi-View Shape Reconstruction,
International Journal of Computer Vision 129 (2021) 284--299.

\bibitem{3Dimmersive1}
D. S. Alexiadis, P. Daras, L. Onural, J. Ostermann, M. A. Magnor, J.-S. Kim, D. Izquierdo,
Reconstruction for 3D immersive virtual environments,
2012 13th International Workshop on Image Analysis for Multimedia Interactive Services, Dublin, Ireland, 2012, pp. 1--4.


\bibitem{MMimage0}
L. Lechelek, S. Horna, R. Zrour, M. Naudin, C. Guillevin, A Hybrid Method for 3D Reconstruction of MR Images,
J. Imaging 2022, 8(4), 103.



\bibitem{JXsurface0}
J. Yang, J. Kim,
Phase-field simulation of multiple fluid vesicles with a consistently energy-stable implicit-explicit method,
Comput. Methods Appl. Mech. Engrg. 417 (2023) 116403.

\bibitem{JXsurface1}
Q. Xia, Y. Liu, J. Kim, Y. Li,
Binary thermal fluids computation over arbitrary surfaces with second-order accuracy and unconditional energy stability based on phase-field model,
J. Comput. Appl. Math. 433 (2023) 115319.

\bibitem{3DpointCloud1}
Y. Lyu, M. Guo, O. Sha, H. Zhang,
3D Point Cloud Surface Reconstruction Based on Divide-and-Conquer Method in Laser Scanner,
2020 J. Phys.: Conf. Ser. 1544 012118.

\bibitem{3DpointCloud2}
Z.M. Bi, L. Wang,
Advances in 3D data acquisition and processing for industrial applications,
Robotics and Computer-Integrated Manufacturing 26(5) (2010) 403-413.

\bibitem{Edelsbrunner1994}
H. Edelsbrunner, E. P. Mücke, Three-dimensional alpha shapes, ACM Trans. Graph. 13 (1994) 43--72.

\bibitem{GCN}
D. Wu, L. Zhou, J. Li, J. Xiong, L. Song,
Explicit 3D reconstruction from images with dynamic graph learning and rendering-guided diffusion, Neurocomputing 601 (2024) 128206.




\bibitem{ACcmame01}
J. Wang, H. Xu, J. Yang, J. Kim,
Fractal feature analysis based on phase transitions of the Allen--Cahn and Cahn--Hilliard equations,
J. Comput. Sci. 72 (2023) 102114.




\bibitem{Osher1988}
S. Osher, J. A. Sethian,
Fronts propagating with curvature-dependent speed: algorithms based on Hamilton--Jacobi formulations,
\textit{J. Comput. Phys.} 79 (1988) 12--49.

\bibitem{Sethian1999}
J. A. Sethian,
\textit{Level Set Methods and Fast Marching Methods},
Cambridge University Press, 1999.

\bibitem{OsherFedkiw2003}
S. Osher, R. Fedkiw,
\textit{Level Set Methods and Dynamic Implicit Surfaces},
Springer, 2003.

\bibitem{Peng1999}
D. Peng, B. Merriman, S. Osher, H. Zhao, M. Kang,
A PDE-based fast local level set method,
\textit{J. Comput. Phys.} 155 (1999) 410--438.

\bibitem{Li2010}
C. Li, C. Xu, C. Gui, M. D. Fox,
Distance Regularized Level Set Evolution and Its Application to Image Segmentation,
\textit{IEEE Trans. Image Process.} 19(12) (2010) 3243--3254.




\bibitem{JSKJWYB02}
J. Wang, H. Xu, W. Jiang, Z. Han, J. Kim,
A novel MF-DFA-Phase-Field hybrid MRIs classification system,
Expert System Appl. 225 (2023) 120071.



\bibitem{JSKJWYB022}
X. Song, Q. Xia, J. Kim, Y. Li,
An unconditional energy stable data assimilation scheme for Navier--Stokes--Cahn--Hilliard equations with local discretized observed data,
Comput. Math. Appl. 164 (2024) 21--33.



\bibitem{JSKJWYB04}
Y. Li, K. Qin, Q. Xia, J. Kim,
A second-order unconditionally stable method for the anisotropic dendritic crystal growth model with an orientation-field,
Appl. Numer. Math. 184 (2023) 512--526.

\bibitem{JSKJWYB041}
W. Xie, Z. Wang, J. Kim, X. Sun, Y. Li,
A novel ensemble Kalman filter-based data assimilation method with an adaptive strategy for dendritic crystal growth,
J. Comput. Phys. 524 (2025) 113711.

\bibitem{face1}
P. Li, W. Li, Y. Tan, H. Fan, Q. Wang,
A phase field fracture model for ultra-thin micro-/nano-films with surface effects,
Int. J. Eng. Sci. 195 (2024) 104004.


\bibitem{face2}
P.-L. Bian, H. Qing, S. Schmauder, T. Yu,
A variationally-consistent phase-field cohensive zone model for mixed-mode fracture with directional energy decomposition scheme and modified-G criterion,
Int. J. Eng. Sci. 210 (2025) 104223.

\bibitem{CCYYY01}
J. Yang,
Unconditionally energy-stable linear convex splitting algorithm for the $L^2$ quasicrystals,
Comput. Phys. Commun. 295 (2024) 108984.


\bibitem{acstruc2}
S. Lai, B. Jiang, Q. Xia, B. Xia, J. Kim, Y. Li,
On the phase-field algorithm for distinguishing connected regions in digital model,
Eng. Anal. Bound. Elem. 168 (2024) 105918.


\bibitem{image001}
J. Wang, Z. Han, J. Kim,
An efficient and explicit local image inpainting method using the Allen--Cahn equation,
Z. Angew. Math. Phys. 75 (2024) 44.


\bibitem{image002}
S. Ham, H. Kim, Y. Hwang, S. Kwak, Jyoti, J. Wang, H. Xu, W. Jiang, J. Kim,
A novel phase-field model for three-dimensional shape transformation,
Comput. Math. Appl. 176 (2024) 67--76.


\bibitem{image003}
H. Kim, S. Kang, G. Lee, S. Yoon, J. Kim,
Shape transformation on curved surfaces using a phase-field model,
Commun. Nonlinear Sci. Numer. Simulat. 133 (2024) 107956.

\bibitem{ACe1}
J. Yang, J. Wang, S. Kwak, S. Ham, J. Kim,
A modified Allen--Cahn equation with a mesh size-dependent interfacial parameter on a triangular mesh,
Comput. Phys. Commun. 304 (2024) 109301.


\bibitem{ACe4}
Z. Liu, X. Li,
A highly efficient and accurate exponential semi-implicit scalar auxiliary variable (ESI-SAV) approach for dissipative system,
J. Comput. Phys. 447 (2021) 110703.

\bibitem{ACe5}
F. Zhang, H.-W. Sun, T. Sun,
Efficient and unconditionally energy stable exponential-SAV schemes for the phase field crystal equation,
Appl. Math. Comput. 470 (2024) 128592.


\bibitem{ACe6}
C. Zhang, J. Ouyang, C. Wang, S.M. Wise,
Numerical comparison of modified-energy stable SAV-type schemes and classical BDF methods on benchmark problems for the functionalized
Cahn--Hilliard equation,
J. Comput. Phys. 423 (2020) 109772.


\bibitem{ACe7}
Y. Qian, Y. Huang, Y. Zhang,
Decoupled, linear and positivity-preserving schemes for a modified phase field crystal system incorporating long-range interactions,
Appl. Math. Comput. 488 (2025) 129089.


\bibitem{ACe8}
Z. Lv, X. Song, J. Feng, Q. Xia, B. Xia, Y. Li,
Reduced-order prediction model for the Cahn--Hilliard equation based on deep learning,
Eng. Anal. Bound. Elem. 172 (2025) 106118.

\bibitem{ACe9}
S. Ham, J. Choi, S. Kwak, J. Kim,
A structure-preserving explicit numerical scheme for the Allen--Cahn equation with a logarithmic potential,
J. Math. Anal. Appl. 538 (2024) 128425.


\bibitem{CCWW01}
C. Yao, F. Zhang, C. Wang,
A scalar auxiliary variable (SAV) finite element numerical scheme for the Cahn--Hilliard--Hele--Shaw system with dynamic boundary conditions,
J. Comput. Math. 42(2) (2024) 544--569.

\bibitem{CCWW02}
Q. Cheng, C. Wang,
Error estimate of second order accurate scalar auxiliary variable (SAV) scheme for the thin film epitaxial models,
Adv. Appl. Math. Mech. 13 (2021) 1318--1354.

\bibitem{CCWW03}
M. Wang, Q. Huang, C. Wang,
A second order accurate scalar auxiliary variable (SAV) numerical method for the square phase field crystal equation,
J. Sci. Comput. 88(2) (2021) 33.

\bibitem{Cheng2020}
Q. Cheng, C. Liu, J. Shen,
A new Lagrange multiplier approach for gradient flows,
Comput. Methods Appl. Mech. Engrg. 367 (2020) 113070.

\bibitem{HuangSurvey2022}
Z. Huang, Y. Wen, Z. Wang, J. Ren, K. Jia,
Surface Reconstruction from Point Clouds: A Survey and a Benchmark,
arXiv:2205.02413 (2022).

\bibitem{SulzerSurvey2023}
R. Sulzer, L. Landrieu, R. Marlet, B. Vallet,
A Survey and Benchmark of Automatic Surface Reconstruction from Point Clouds,
arXiv:2301.13656 (2023).

\bibitem{VIPSS2019}
Z. Huang, N. Carr, T. Ju,
Variational implicit point set surfaces,
ACM Trans. Graph. 38(4) (2019) 124:1--124:13.

\bibitem{iPSR2022}
F. Hou, C. Wang, W. Wang, H. Qin,
Iterative Poisson Surface Reconstruction (iPSR) for Unoriented Points,
ACM Trans. Graph. 41(4) (2022) 128:1--128:13.

\bibitem{DeepSDF2019}
J.J. Park, P. Florence, J. Straub, R. Newcombe, S. Lovegrove,
DeepSDF: Learning Continuous Signed Distance Functions for Shape Representation,
Proc. IEEE/CVF CVPR (2019).

\bibitem{OccNet2019}
L. Mescheder, M. Oechsle, M. Niemeyer, S. Nowozin, A. Geiger,
Occupancy Networks: Learning 3D Reconstruction in Function Space,
Proc. IEEE/CVF CVPR (2019).

\bibitem{SAL2020}
M. Atzmon, Y. Lipman,
SAL: Sign Agnostic Learning of Shapes from Raw Data,
Proc. IEEE/CVF CVPR (2020).

\bibitem{IGR2020}
A. Gropp, L. Yariv, N. Haim, M. Atzmon, Y. Lipman,
Implicit Geometric Regularization for Learning Shapes,
Proc. ICML, PMLR 119 (2020) 3789--3799.

\bibitem{NeuralPull2021}
B. Ma \textit{et al.},
Neural-Pull: Learning Signed Distance Functions from Point Clouds by Learning to Pull Space onto Surfaces,
Proc. ICML, PMLR 139 (2021).

\bibitem{Points2Surf2020}
P. Erler, P. Guerrero, S. Ohrhallinger, M. Wimmer, N.J. Mitra,
Points2Surf: Learning Implicit Surfaces from Point Cloud Patches,
Proc. ECCV (2020).


\bibitem{allencahn}
S.M. Allen, J.W. Cahn,
Ground state structures in ordered binary alloys with second neighbor interactions,
Acta Metall. 20 (1972) 423--433.


\bibitem{JJYangJCP}
Y. Zhao, D. Cai, J. Yang,
Second-order accurate and unconditionally stable algorithm with unique solvability for a phase-field model of 3D volume reconstruction,
J. Comput. Phys. 504 (2024) 112873.

\bibitem{Cai2025}
D. Cai, B. Fu, R. Gao, X. Kong, J. Yang,
Phase-field computation for 3D shell reconstruction with an energy-stable and uniquely solvable BDF2 method,
\textit{Comput. Math. Appl.} 189 (2025) 1--23.

\bibitem{Kong2025}
X. Kong, R. Gao, B. Fu, D. Cai, J. Yang,
Two lower boundedness-preservity auxiliary variable methods for a phase-field model of 3D narrow volume reconstruction,
\textit{Commun. Nonlinear Sci. Numer. Simul.} 143 (2025) 108649.



\bibitem{c3podata}
Gambody, C-3PO 3D Printing Model. Available at: https://www.gambody.com/premium/c-3po.

\bibitem{c3pomovie}
Wikipedia, C-3PO (droid) image. Available at: https://upload.wikimedia.org/wikipedia/en/5/5c/C-3PO\_droid.png.



\bibitem{Vaderdata1}
Sketchfab, Darth Vader Helmet (3D model). Available at: https://sketchfab.com/3d-models/darth-vader-helmet-336c207bda9b4528b16e1dfa61007986.



\bibitem{SAVV1}
J. Shen, J. Xu, J. Yang,
The scalar auxiliary variable (SAV) approach for gradient flows,
J. Comput. Phys. 353 (2018) 407--416.


\bibitem{SAVV2}
J. Yang, J. Kim,
Consistently and unconditionally energy-stable linear method for the diffuse-interface model of narrow volume reconstruction,
Eng. Comput. 40 (2024) 2617--2627.


\end{thebibliography}
\end{document}